\documentclass[11pt,reqno,a4wide]{article}
\usepackage{a4wide}
\usepackage{eurosym}
\usepackage{comment}
\usepackage{mathptmx}
\usepackage{amssymb}
\usepackage{amsfonts}
\usepackage{amsthm}
\usepackage{amsmath}
\usepackage{graphicx}
\usepackage{caption}
\usepackage{shadow}
\usepackage[all]{xy}
\usepackage{enumerate}
\usepackage{makeidx}
\usepackage{mathrsfs}
\usepackage{upgreek}
\usepackage{stmaryrd}
\usepackage[utf8]{inputenc}
\usepackage{array}
\usepackage{tikz-cd}
\usepackage{thumbpdf}
\usepackage[pagebackref]{hyperref}
\usepackage[notref,notcite]{showkeys}
\usepackage{textcomp}
\usepackage{chapterbib}
\usepackage{graphics}
\usepackage{longtable}
\usepackage{bm}
\usepackage{adjustbox}
\usepackage{centernot}
\usepackage{setspace}
\usepackage{imakeidx}
\usepackage{tikz}
\usepackage{tikz-cd}
\usepackage{float}

\usetikzlibrary{positioning}
\makeindex
\newtheorem{theorem}{Theorem}[section]
\newtheorem{lem}[theorem]{Lemma}
\newtheorem{thm}[theorem]{Theorem}
\newtheorem{prop}[theorem]{Proposition}
\newtheorem{cor}[theorem]{Corollary}
\theoremstyle{definition}
\newtheorem*{Beweis}{Proof}

\newtheorem{defn}[theorem]{Definition}

\newtheorem{rem}[theorem]{Remark}
\newtheorem{rems}[theorem]{Remarks}
\newtheorem{punto}[theorem]{}
\theoremstyle{remark}

\newtheorem{ex}[theorem]{Example}
\newtheorem{exs}[theorem]{Examples}
\input{tcilatex}
\begin{document}

\title{{Semirings Induced by Semilattices}}
\author{Jawad Abuhlail \\
\text{abuhlail@kfupm.edu.sa; abuhlail@yahoo.de}\\
Department of Mathematics\\
Center for Intelligent Secure Systems\\
King Fahd University of Petroleum $\&$ Minerals\\
31261 Dhahran, KSA}
\date{\today}
\maketitle

\begin{abstract}
We investigate \emph{semirings} induced by join-semilattice
(meet-semilattices) with a degenerate multiplication on all the subset that
excludes the multiplicative unity. The way those semirings are defined
enables a big impact of the semilattices we start with on the associated
semirings: many of their order-theoretic properties reflect in algebraic
properties of the induced semirings. This enabled constructing examples of
semirings with desired properties as well as providing counterexamples to
wrong claims in old and recent papers/preprints. In particular, we construct
families of examples of semirings induced by \emph{non-modular}
join-semilattices for which the lattice of subtractive ideals is \emph{%
non-modular} disproving a claim of Henriksen who introduced this important
class of ideals of semirings.
\end{abstract}

\section*{Introduction}

\bigskip

A semiring is (roughly)\ an associative ring \emph{not necessarily with
subtraction}. The earliest roots of semirings in modern algebra go back to
Abel 1826 and Dedekind 1894 (cf. \cite{Gla2002}). The term \emph{semiring}
was formally introduced by H. S. Vandiver in 1934 \cite{Van1934} as an
algebraic system that generalizes both associative rings and distributive
bounded lattices.

\bigskip

Several types of semirings are related, in a sense or another, to
(semi)lattices in a way that their fundamental \emph{algebraic }properties
are reflection of \emph{order theoretic} properties of the probably hidden
(semi)lattice structures that induced them. A \emph{nontrivial }bounded
lattice $(L,\wedge ,\vee ,0,1)$ is a semiring with addition/multiplication
given by $\vee /\wedge $ (or $\wedge /\vee $) if and only if $L$ is
distributive (cf. Theorem \ref{classical-lat}). Every \emph{additively
idempotent} semiring $S$ has a structure of a join-semilattice induced by
the natural order $x\leq _{+}y\Longleftrightarrow x+y=y.$ A finite \emph{%
proper} hemiring $H$ is \emph{congruence-simple} if and only if $\left\vert
H\right\vert \leq 2$ or $H$ is isomorphic to a dense subhemiring of the
endomorphism semiring $End(M)$ of an additively idempotent monoid $M$ \cite[%
Theorem 1.7]{Zum2008}.

\bigskip

In this paper, we associate to each join-semilattice $(L,\vee )$ with $%
\left\vert L\right\vert \geq 2$ and fixed $c\in L,$ an additively idempotent
demiring $D_{\vee }^{c}(L)$ in a way that the lattice $\func{Idl}(L)$ of 
\emph{order ideals} of $U_{c}:=[c\uparrow )$ and the lattice $\mathfrak{K}%
(D_{\vee }^{c}(L))$ of subtractive ideals of the associated demiring are
isomorphic. If $(L,\vee ,0)$ is bounded, then $H_{\vee }(L):=D_{\vee
}^{0}(L) $ is a hemiring. Moreover, one can adjoin to $L$ a \emph{%
join-irreducible} top element $1^{\sharp }\notin L$ (if it does not have
such top) and construct a Krull zero-dimensional local semiring $S_{\vee
}(L^{\sharp })$ with underlying set $L^{\sharp }:=L\cup \{1^{\sharp }\}$ in
a way that the lattice $\func{Idl}(L^{\sharp })=\func{Idl}(L)\cup
\{L^{\sharp }\}$ of \emph{order ideals} of $L^{\sharp }$ and the lattice $%
\mathfrak{K}(S_{\vee }(L^{\sharp }))$ of subtractive ideals of the induced
semiring are isomorphic. One big advantage is that these constructions can
be dualized to induce a structure of a demiring (hemiring) to each (bounded)
meet-semilattice $(M,\wedge )$ with a fixed $c\in M$ ($c=\top $, the top
element of $M$), as well as a semiring in case $(M,\wedge ,\top )$ is
bounded with a meet-irreducible bottom $\bot \neq \top .$

\bigskip

Since the multiplication of a demiring (hemiring) induced by a (bounded)\
join-semilattice is constant, and the multiplication of a semiring induced
by a bounded join-semilattice with a join-irreducible top is constant \emph{%
except} when it involves the top element, we have a natural \emph{machinery}
of constructing examples/counterexamples of additively idempotent
commutative demirings and hemirings (semirings) having/missing a specific
property: we impose conditions on an appropriately chosen join-semilattice
whose \emph{order} properties reflect in the desired/forbidden \emph{%
algebraic} properties of the induced semiring. Using this machinery, and in
addition to numerous examples illustrating our results, we managed to
provide \emph{families} of semirings induced by semilattices which serve as 
\textbf{counterexamples} to two \emph{wrong claims} that appeared recently
in \cite{GD2024}:\ the modularity of the lattice of subtractive ideals of 
\emph{any} commutative semiring \cite[Proposition 2.6]{GD2024}, and that
commutativity of the subtractive closure of ideals of \emph{any} commutative
semiring with intersections \cite[Lemma 2.1 (7)]{GD2024}.

\bigskip

In Section 1, we collect some preliminaries from Semiring Theory and Lattice
Theory with the convention that the underlying set of each join-semilattice,
lattice, demiring, hemiring and semiring \emph{we start with} has at least
two elements. In Section 2, we study demirings, hemirings and semirings
induced by join-semilattices (meet-semilattices). Let $(L,\vee )\ $be a
join-semilattice. In Proposition \ref{semiring-S0}, we construct a demiring $%
D_{\vee }^{c}(L)$ associated to any fixed $c\in L$ with constant
multiplication, and show that $H_{\vee }(L)=D_{\vee }^{0}(L)$ is a hemiring
with zero-multiplication if and only if $L$ is bounded with bottom $\bot $
and $\bot =0;$ and if $(L,\vee ,0)$ is bounded with a \emph{join-irreducible}
top element $1\neq 0,$ then we obtain a structure of semiring $S_{\vee }(L)$
with zero-multiplication on $L\backslash \{1\}.$ Proposition \ref%
{j-semiring-local} proves that $S_{\vee }(L)$ is local and that its unique
maximal ideal $\mathfrak{m}=L\backslash \{1\}\ $is subtractive. Proposition %
\ref{ideals-join} and Corollary \ref{ideals-join-S} are the central
technical results in this paper as they provide us with exact descriptions
of the ideals and subtractive ideals of the demiring $D_{\vee }^{c}(L)$ (the
semiring $S_{\vee }(L)$). In Theorem \ref{mod-S-join} (Theorem \ref%
{dist-S-join}), we prove that $\mathfrak{K}(D_{\vee }^{c}(L))$ is modular
(distributive) if and only if the join-semilattice $L$ is modular
(distributive). This result is extended in Theorem \ref{1-hash} to the
semiring $\mathfrak{K}(S_{\vee }(L))$ is case $(L,\vee ,0)$ is a bounded
join-semilattice with join-irreducible $1\neq 0.$ These basic results are
then dualized in Proposition \ref{semiring-meet} to the demiring $D_{\wedge
}^{c}(M)$ induced by a meet-semilattice $(M,\wedge )$ and fixed $c\in M,$
the hemiring $H_{\wedge }(M)$ (the semiring $S_{\wedge }(M)$) induced by a
bounded meet-semilattice $(M,\wedge ,\top )\ $with a meet-irreducible bottom
element $\bot \neq \top $). Then we investigate several types of \emph{%
simplicity} for demirings, hemirings and demirings induced by
join-semilattices. Let $(L,\vee )\ $be a join-semilattice. Theorem \ref%
{join-simplicity} provides separate characterizations for the demiring $%
D_{\vee }^{c}(L)$ to be congruence-simple (resp. ideal-simple, $k$%
-ideal-simple). If $(L,\vee ,0)$ is bounded (with a join-irreducible $1\neq
0 $), then Theorem \ref{H-CC-primes} (Theorem \ref{join-c-simple}), shows
that these three types of simplicity are equivalent to $H_{\vee }(L)\cong 
\mathfrak{B,}$ the \emph{Boolean hemiring }zero-multiplication ($S_{\vee
}(L)\cong \mathbb{B}$ the Boolean semiring).

\bigskip

In Section 3, we study what we call Kuratowski demirings and nucleus\
demiring whose definitions are inspired by \cite{Nas}. We call a demiring $%
(D,+,\cdot )$ a \emph{Kuratowski demiring} (resp. a \emph{nucleus demiring}%
)\ iff the subtractive closure of its ideals commutes with binary sums
(binary intersections). Given a join-semilattice $(L,\vee ,0)$ (with a
join-irreducible $1\neq 0$), we show in Theorem \ref{j-K-dist} that $H_{\vee
}(L)$ is a Kuratowski hemiring ($S_{\vee }(L)$) if and only if $L$ is
distributive. It shows also that if $(L,\vee )$ is a distributive
join-semilattice and $c\in L$ is fixed, then the distributivity of $L$
implies that $D_{\vee }^{c}(L)$ is a Kuratowski demiring (the converse is
not true in general as illustrated by Example \ref{N5-not-Kurat}). On the
other hand, Theorem \ref{nuc-dem} (Theorem \ref{j-Nuc-dist}) shows that the
demiring $D_{\vee }^{c}(L)$ (the semiring $S_{\vee }(L)$) is nucleus if and
only if $\left\vert L\right\vert =2$ ($\left\vert L\right\vert =2$ or $%
\left\vert L\right\vert =3$). While counterexamples to \cite[Lemma 2.1 (7)]%
{GD2024} already exist (cf. \cite{SGBS} and \cite{Nas}), our theorem does
not only give infinitely many \textbf{counterexamples} to \cite[Lemma 2.1 (7)%
]{GD2024}, but also shows when that claim holds for our structures.

\bigskip

Section 4 is devoted to chain conditions for semirings induced by
semilattices. Theorem \ref{j-semiring-NA} characterizes the demirings $%
D_{\vee }^{c}(L)$ that are Noetherian (Artinian)\ as those for which the
class $\mathfrak{I}_{\vee }^{c}(U_{c})$ (of all subsets of $L$ containing $c$
and closed under the join) satisfies the ACC (DCC), and demirings which are $%
k$-Noetherian ($k$-Artinian)\ as those for which the class $\func{Idl}%
(U_{c}) $ of \emph{order ideals} of $U_{c}$ satisfies the ACC (DCC). These
characterizations extend to the hemiring $H_{\vee }(L)=D_{\vee }^{0}(L)$
induced by bounded join-semilattice $(L,\vee ,0)\ $and to the semiring $%
S_{\vee }(L)$ induced by a bounded join-semilattice with join-irreducible
top $1\neq 0.$ In the special case $L=C$, a chain, we demonstrate in
Corollary \ref{chain-4-types} that $S_{\vee }(C)$ is Noetherian/Artinian if
and only if $C$ is finite. Moreover, we prove that $S_{\vee }(C)$ is $k$%
-Noetherian ($k$-Artinian)\ if and only if $C$ satisfies the ACC (DCC). This
enables us to provide examples of a \emph{Krull zero-dimensional}
commutative semirings which are $k$-Noetherian but not $k$-Artinian, or $k$%
-Artinian but not $k$-Noetherian. This illustrates that \emph{proper}
commutative semirings are quite weird and behave quite differently from
commutative rings with respect to chain conditions.

\section{Preliminaries}

\bigskip

In this section, we present the main definitions from the Theory of
Semirings and Semimodules (following Golan \cite{Gol1999})\ and Lattice
Theory (following Gr\"{a}zer \cite{Gra2010}).

\bigskip

There is no unified definition of a semiring in the literature. To avoid
confusion, we give each of the main types of \emph{semirings} in the
literature a special name to distinguish it from the others and to make any 
\emph{upgrade} clear.

\begin{punto}
A \textbf{demiring} is a datum $(S,+,\cdot ),\ $where $(S,+)$ is a
commutative semigroup, $(S,\cdot )\ $is a semigroup, and multiplication
distributes over addition from both sides. A \textbf{pre-hemiring} is a
demiring $(S;+,0;\cdot )$ in which $(S,+,0)$ is a monoid. A \textbf{hemiring}
is a pre-hemiring in which $0$ is \emph{absorbing} in the sense that $s\cdot
0=0=0\cdot s$ for all $s\in S.$ A \textbf{semiring} is a hemiring $%
(S;+,0;\cdot ,1)\ $in which $(S,\cdot ,1)$ is a monoid and $1\neq 0.$ To
avoid the trivial case, we assume $\left\vert S\right\vert \geq 2$ for any
demiring in this paper. A demiring $(S,+,\cdot )$ is said to be \textbf{%
commutative} iff $(S,\cdot )$ is commutative, \textbf{additively idempotent}
iff $s+s=s$ for every $s\in S,$ \textbf{multiplicatively idempotent} iff $%
s\cdot s=s$ for every $s\in S,$ and \textbf{idempotent} iff $S$ is
additively idempotent and multiplicatively idempotent. An ideal of a
demiring $(S,+,\cdot )$ is a \emph{non-empty} subset $I\subseteq S$ such
that $I+I\subseteq I,$ $SI\subseteq I$ and $IS\subseteq I.$ With $\mathfrak{I%
}(S)$ we denote the class of ideals of $S.$ A \textbf{proper demiring} is a
demiring which is \emph{not }a ring.
\end{punto}

\bigskip

\begin{punto}
Let $(S,+,\cdot )\ $be a demiring. For any subset $A\subseteq S,$ we define
the \textbf{subtractive closure} of $A$ as%
\begin{equation*}
\overline{A}:=\{s\in S\text{ }|\text{ }s+a\in A\text{ for some }a\in A\}.
\end{equation*}%
We say that $A$ is \textbf{subtractive} iff $\overline{A}\subseteq A$ ($%
\overline{A}=A$ if $S$ is pre-hemiring). The subtractive closure of any
ideal is an ideal. If $I\subseteq S$ is an ideal, then $I$ is subtractive if
and only if $\overline{I}=I.$ A subtractive ideal is also called a $k$%
\textbf{-ideal, }and the class of class of subtractive ideals ($k$-ideals)
of $S$ is denote by $\mathfrak{K}(S).$ A commutative hemiring is said to be 
\textbf{subtractive} iff every ideal of $S$ is subtractive.
\end{punto}

\bigskip

\begin{punto}
Let $(S,+,\cdot )$ be a demiring. For a \emph{non-empty} subset $\varnothing
\neq A\subseteq S,$ define the \textbf{additive closure} of $A$ as%
\begin{equation*}
A^{+}:=\{x+s\mid x\in A\text{ and }s\in S\}.
\end{equation*}%
We say that $A$ is \textbf{additive} iff $A^{+}\subseteq A$ ($A=A^{+}$ if $S$
is a pre-hemiring). With $\mathcal{A}(S)$ we denote the collection of all
additive subsets of $S.$ A \textbf{coideal} of $S$ is an additive subset $%
C\subseteq S$ that is closed under multiplication. With $\mathcal{C}(S)$ we
denote the collection of all coideals of $S.$ A \textbf{bi-ideal} of $S$ is
an ideal that is also a coideal. With $\mathcal{B}(S)=\mathfrak{I}(S)\cap 
\mathcal{C}(S)$ we denote the collection of all bi-ideals of $S.$
\end{punto}

\begin{defn}
We call a demiring $(S,+,\cdot )$ with $\left\vert S\right\vert \geq 2:$

\textbf{congruence-simple} iff the class of all congruences of $S$ is $\func{%
Cong}(S)=\{\Delta ,\nabla \},$ where $\Delta $ is the diagonal (called the
trivial congruence) and $\nabla =S\times S$ (called the universal
congruence);

\textbf{ideal-simple} iff the only ideal $I$ of $S$ with $\left\vert
I\right\vert \geq 2$ is $S;$

$k$\textbf{-ideal-simple} iff the only $k$-ideal $I$ of $S$ with $\left\vert
I\right\vert \geq 2$ is $S.$
\end{defn}

\begin{punto}
Let $(S,+,\cdot )\ $be a hemiring. The class $(\mathfrak{I}(S),\cap ,+)$ of 
\emph{all} ideals of $S$ is a complete bounded lattice with bottom element $%
\{0\}\ $and top element $S.$ While $(\mathfrak{K}(S),\cap )\subseteq (%
\mathfrak{I}(S),\cap )$ is a meet-subsemilattice, $\mathfrak{K}(S)$ is \emph{%
not} even closed under addition (cf. Example \ref{N5+1-no-KN}). To overcome
this obstacle, one adjusts the join on $\mathfrak{K}(S)$ and defines $I\vee
_{\mathfrak{K}}J:=\overline{I+J}$ (and $\wedge _{\mathfrak{K}}:=\cap $). One
can check easily that $(\mathfrak{K}(S),\wedge _{\mathfrak{K}},\vee _{%
\mathfrak{K}})$ is a lattice.
\end{punto}

\begin{punto}
We call a partially ordered set $(L,\leq )$

\textbf{join-semilattice (upper semilattice)\ }iff any two-element set $%
\{a,b\}$ has a \emph{least upper bound} $a\vee b:=\sup \{a,b\}$ called the 
\textbf{join};

\textbf{meet-semilattice (lower semilattice) }iff any two-element set $%
\{a,b\}$ has a \emph{greatest lower bound} $a\wedge b:=\inf \{a,b\}$ called
the \textbf{meet}.

We call a \emph{join-semilattice} (meet-semilattice)\ 

\textbf{bounded} iff it has a \emph{bottom} $0$ (\emph{top} $\top $);

\textbf{complete} iff arbitrary joins (meets) exist.
\end{punto}

\begin{rem}
Let $(L,\vee )$ ($(L,\wedge )$) be a join-semilattice (meet-semilattice).
Arbitrary joins (meets)\ include the empty ones, whence a complete
join-semilattice (meet-semilattice)\ has bottom $0:=\bigvee \varnothing $ ($%
\top :=\bigwedge \varnothing $) and top $\bigvee L$ (bottom $\bot
:=\bigwedge L$).
\end{rem}

\begin{punto}
A \textbf{lattice} is a poset $(L,\leq )$ in which every two-element set $%
\{a,b\}$ has a meet and a join. A lattice $L$ is bounded (complete) iff both
semilattices $(L,\wedge )\ $and $(L,\vee )$ are bounded (complete). Every
complete lattice is clearly bounded.
\end{punto}

\begin{defn}
Let $\mathcal{L}=(L,\wedge ,\vee ,0,1)$ be a complete lattice. An element $%
c\in L$ is said to be \textbf{compact} iff for every subset $A\subseteq L:$%
\begin{equation*}
c\leq \bigvee A\quad \Longrightarrow \quad c\leq \bigvee F\text{ for some }F%
\underset{\text{finite}}{\subseteq }A;
\end{equation*}%
We say that $\mathcal{L}$ is \textbf{algebraic} iff every element of $L$ is
a join of a (possibly infinite) set of compact elements.
\end{defn}

\begin{punto}
Let $(P,\leq )$ b a poset. For any non-empty subset $A\subseteq P,$ we set%
\begin{eqnarray*}
\lbrack A &\uparrow &):=\{x\in P\text{ }|\text{ }a\leq x\text{ for some }%
a\in A\};\text{ } \\
( &\downarrow &A]:=\{x\in P\text{ }|\text{ }x\leq a\text{ for some }a\in A\}.
\end{eqnarray*}%
If $A=\{a\},$ a singleton, then we set 
\begin{equation*}
U_{a}:=[a\uparrow ):=[\{a\}\uparrow )\text{ and }D_{a}:=(\downarrow
a]:=(\downarrow \{a\}].
\end{equation*}%
We say that a \emph{non-empty} subset $\varnothing \neq A\subseteq P$ is

\textbf{up-set} iff $[A\uparrow )=A,$ equivalently $[a\uparrow )\subseteq A$
for every $a\in A;$

\textbf{down-set} iff $(\downarrow A]=A,$ equivalently $(\downarrow
a]\subseteq A$ for every $a\in A.$

With $\mathfrak{U}(P)$ ($\mathfrak{D}(P)$) we denote the collection of all
up-sets (down-sets)\ of $P.$
\end{punto}

\begin{punto}
Let $(L,\vee )$ be a join-semilattice. A \emph{non-empty} subset $D\subseteq
L$ is called an \textbf{order ideal} iff $D$ is a down-set that is closed
under binary joins, equivalently iff%
\begin{equation}
\forall \text{ }a,b\in L:\text{ }a\vee b\in D\Longleftrightarrow a,b\in D.
\label{Gra-ideal}
\end{equation}%
With $\func{Idl}(L)$ we denote the collection all order ideals of $L.$ For
any $d\in D,\ $the down-set $(\downarrow d]$ is clearly an order ideal,
called the \textbf{principal order ideal} generated by $d.$
\end{punto}

\begin{lem}
\label{join-sem-idl}Let $(L,\vee )$ be a join-semilattice.

\begin{enumerate}
\item $\func{Idl}(L)$ is a join-semilattice with join%
\begin{equation}
I\vee _{\func{Idl}(L)}J:=\{x\in L\text{ }\mid \text{ }x\leq a\vee b\text{
for some }a\in I,b\in J\}.  \label{j-idl}
\end{equation}%
Moreover, $\func{Idl}(L)$ is closed under arbitrary \emph{non-empty} joins.

\item If $L$ is bounded, then $\func{Idl}(L)$ is a \emph{complete algebraic}
lattice with meet the intersection and whose set of compact elements is
precisely the set of principal order ideals $\{(\downarrow d]$ $|$ $d\in
L\}. $
\end{enumerate}
\end{lem}

\begin{Beweis}
Let $(L,\vee )\ $be a join-semilattice.

\begin{enumerate}
\item For any \emph{non-empty} collection $(A_{\lambda })_{\Lambda
}\subseteq \func{Idl}(L),$ consider%
\begin{equation}
\bigvee\limits_{\lambda \in \Lambda }A_{\lambda }:=\{x\in L\text{ }|\text{ }%
x\leq a_{1}\vee \cdots \vee a_{n},\text{ }n\geq 1,\text{ }a_{1},\cdots
,a_{n}\in \bigcup\limits_{\lambda \in \Lambda }A_{\lambda }\}.
\label{join-idl}
\end{equation}%
It's obvious that $\bigvee\limits_{\lambda \in \Lambda }A_{\lambda }$
satisfies (\ref{Gra-ideal}), whence an order ideal of $L.$ So, $\func{Idl}%
(L) $ is closed under arbitrary \emph{non-empty} joins.

\item This is \cite[Theorem 42]{Gra2010}\emph{.}$\blacksquare $
\end{enumerate}
\end{Beweis}

The following example shows that the assumption that the join-semilattice in
Lemma \ref{join-sem-idl} is bounded \emph{cannot} be dropped.

\begin{ex}
Consider the join-semilattice%
\begin{equation*}
L:a,b<t,\quad a\parallel b.
\end{equation*}%
Notice that $L$ is not bounded. Moreover, $\{a\},$ $\{b\}\in \func{Idl}(L)$
but $\{a\}\cap \{b\}=\varnothing \notin \func{Idl}(L).$ So, $\func{Idl}(L)$
is \emph{not} even a lattice.$\blacksquare $
\end{ex}

\begin{punto}
Let $(M,\wedge ;\top )$ be a meet-semilattice with a greatest element $\top
. $ A subset $U\subseteq M$ is called a \textbf{filter} iff $U$ is a \emph{%
non-empty} up-set that is closed under binary meets, equivalently%
\begin{equation}
\forall \text{ }a,b\in M:\text{ }a\wedge b\in U\Longleftrightarrow a,b\in U.
\label{Gra-filter}
\end{equation}%
With $\func{Fil}(M)=\func{Idl}(M^{\mathrm{op}})$ we denote the collection
all filters of $M.$ For every $c\in M,$ the up-set $U_{c}:=[c\uparrow )$ is
clearly a filter, called the \textbf{principal filter} generated by $c.$
\end{punto}

\begin{defn}
(cf. \cite{Gra2010}, \cite{Rho1975}) Let $(L,\leq )$ be a poset.

\begin{enumerate}
\item If $(L,\vee )$ is a join-semilattice, then we call $L:$

\textbf{modular join-semilattice} iff for all $u,v,w\in L$ with $u\leq w\leq
u\vee v$ there exists $y\in L$ with $y\leq v$ and $w=u\vee y;$

\textbf{distributive join-semilattice} iff for all $u,v,w\in L$ with $u\leq
v\vee w$ there exist $y,z\in L$ with $y\leq v,$ $z\leq w$ and $u=y\vee z.$

\item If $(L,\wedge )$ is a meet-semilattice, then we call $L:$

\textbf{modular meet-semilattice} iff for all $u,v,w\in L$ with $u\wedge
v\leq w\leq u$ there exists $y\in L$ with $v\leq y$ and $w=u\wedge y;$

\textbf{distributive meet-semilattice} iff for all $u,v,w\in L$ with $%
v\wedge w\leq u$ there exist $y,z\in L$ with $v\leq y,$ $w\leq z$ and $%
u=y\wedge z.$

\item If $(L,\wedge ,\vee )$ is a lattice, then we call $L:$

\textbf{modular lattice} iff for all $A,B,C\in L$ with $A\leq C$ we have%
\begin{equation*}
(A\vee B)\wedge C=A\vee (B\wedge C);
\end{equation*}

\textbf{distributive lattice} iff for all $A,B,C\in L$ we have%
\begin{eqnarray*}
(A\vee B)\wedge C &=&(A\wedge C)\vee (B\wedge C);\text{ } \\
\text{(equivalently }(A\wedge B)\vee C &=&(A\vee C)\wedge (B\vee C)\text{).}
\end{eqnarray*}
\end{enumerate}
\end{defn}

The following result is well known (cf. \cite[Corollary 1.2]{Cor1974}):

\begin{lem}
\label{equiv-lat-semi}Let $(L,\wedge ,\vee )$ be a lattice.

\begin{enumerate}
\item $L$ is modular as a lattice $\Longleftrightarrow $ $(L,\vee )$ is
modular as a join-semilattice $\Longleftrightarrow $ $(L,\wedge )$ is
modular as a meet-semilattice.

\item $L$ is distributive as a lattice $\Longleftrightarrow $ $(L,\vee )$ is
distributive as a join-semilattice $\Longleftrightarrow $ $(L,\wedge )$ is
distributive as a meet-semilattice.
\end{enumerate}
\end{lem}

\begin{punto}
In both diagrams below, the underlying set is $\{e,a,b,c,t\},$ where $e$ is
the least element and $t$ is the greatest element. The diagrams are Hasse
diagrams: an edge joins two elements precisely when one covers the other.
The covering relations are 
\begin{eqnarray*}
N_{5} &:&e<a<c<t;\qquad e<b<t;\qquad b\parallel a,c; \\
M_{3} &:&e<a,b,c;\qquad a,b,c<t;\qquad a\parallel b\parallel c.
\end{eqnarray*}%
\begin{figure}[th]
\centering
\begin{minipage}{0.47\textwidth}
\centering
\begin{tikzpicture}[scale=1.0, every node/.style={circle,draw,fill=white,inner sep=1.5pt,minimum size=7mm}]
  \node (e) at (0,0) {$e$};
  \node (a) at (-1.15,1) {$a$};
  \node (c) at (-1.15,2) {$c$};
  \node (t) at (0,3) {$t$};
  \node (b) at (1.15,1.5) {$b$};
  \draw (e)--(a)--(c)--(t);
  \draw (e)--(b)--(t);
\end{tikzpicture}
\label{pentagon}
\caption{The \textbf{pentagon} lattice $N_5$.}
\end{minipage}
\hfill 
\begin{minipage}{0.47\textwidth}
\centering
\begin{tikzpicture}[scale=1.0, every node/.style={circle,draw,fill=white,inner sep=1.5pt,minimum size=7mm}]
  \node (e) at (0,0) {$e$};
  \node (a) at (-1.15,1.35) {$a$};
  \node (b) at (0,1.65) {$b$};
  \node (c) at (1.15,1.35) {$c$};
  \node (t) at (0,3) {$t$};
  \draw (e)--(a);
  \draw (e)--(b);
  \draw (e)--(c);
  \draw (a)--(t);
  \draw (b)--(t);
  \draw (c)--(t);
\end{tikzpicture}
\label{diamond}
\caption{The \textbf{diamond} lattice $M_3$.}
\end{minipage}
\label{pent-diamond}
\end{figure}
\end{punto}

\begin{lem}
\label{mod-dist}Let $L$ be \emph{any} lattice.

\begin{enumerate}
\item $L$ is modular if and only if $L$ contains \emph{no sublattice}
isomorphic to $N_{5}.$

\item $L$ is distributive if and only if $L$ contains \emph{no sublattice}
isomorphic to $N_{5}$ or $M_{3}.$
\end{enumerate}
\end{lem}

The following result are \emph{folklore} (cf. \cite{Gol1999}).

\begin{thm}
\label{classical-lat}The following are equivalent for a \emph{non-trivial}
lattice $(L,\wedge ,\vee ):$

\begin{enumerate}
\item $L$ is distributive;

\item $Lat_{\vee }(L):=(L,\vee ,\wedge )$ is an (idempotent) demiring with
addition $\vee $ and multiplication $\wedge ;$

\item $Lat_{\wedge }(L):=(L,\wedge ,\vee )$ is an (idempotent) demiring with
addition $\wedge $ and multiplication $\vee .$

Moreover, if the lattice $L$ is bounded, with bottom $0$ and top $1,$ and
distributive, then we have an isomorphism of semirings:%
\begin{equation*}
Lat_{\vee }(L):=(L;\vee ,0;\wedge ,1)\cong Lat_{\wedge }(L):=(L;\wedge
,1;\vee ,0).
\end{equation*}
\end{enumerate}
\end{thm}

\section{Semirings Induced by Join-Semilattices}

Every algebra in this paper is nontrivial in the sense that it has at least
two elements.

\begin{defn}
An element $z$ in a bounded join-semilattice $(L,\vee ;\bot )$
(meet-semilattice $(M,\wedge ,\top )$) is said to be join-irreducible\textbf{%
\ (meet-irreducible}) \cite[p. 102]{Gra2010} iff%
\begin{eqnarray}
z &\neq &\bot \text{ and for any }x,y\in L\text{ with\ }x\vee y=z,\text{ we
have }x=z\text{ or }y=z;  \label{top-join-prime} \\
z &\neq &\top \text{ and for any }x,y\in M\text{ with\ }x\wedge y=z,\text{
we have }x=z\text{ or }y=z.  \label{bot-meet-pfime}
\end{eqnarray}
\end{defn}

\begin{lem}
\label{cover}Let $(L,\wedge ,\vee ,0,1)$ be a \emph{finite} bounded lattice
with $1\neq 0.$

\begin{enumerate}
\item $1$ is join-irreducible $\Longleftrightarrow $ $1$ has a unique lower
cover.

\item $0$ is meet-irreducible $\Longleftrightarrow $ $0$ has a unique upper
cover.
\end{enumerate}
\end{lem}

The following example shows that the assumption that $L$ is finite in Lemma %
\ref{cover} cannot be dropped.

\begin{ex}
Consider the real interval $[0,1]$ as a lattice ordered in the usual way.
Since $\vee _{L}=\max $ and $\wedge _{L}=\min ,$ we know that $1$ is
join-irreducible and $0$ is meet-irreducible. However, \emph{neither} $1$
has a lower cover \emph{nor} $0$ has an upper cover.

On the other hand, Consider the bounded lattice $(L,\vee ,\wedge ,0,1)$
where 
\begin{eqnarray*}
L &=&\{0,1,c\}\cup \{a_{1},a_{2},\cdots \}\cup \{b_{1},b_{2},\cdots \},\text{
}0<x<1\text{ for all }x\in L, \\
\text{ }c &\parallel &a_{n}\parallel b_{m}\text{ for all }n,m,\qquad
a_{n}<a_{n+1},\text{ }b_{n}<b_{n+1}\text{ for all }n.
\end{eqnarray*}%
Although $c$ is the unique lower cover of $1$ in $L,$ it follows from $%
a_{1}\vee b_{1}=1$ that $1$ is \emph{not} join-irreducible. In the dual
lattice $L^{\mathrm{op}},$ $c$ is a unique upper cover of $0_{L^{\mathrm{op}%
}}=1$ but $0_{L^{\mathrm{op}}}$ is \emph{not} meet-irreducible.$\blacksquare 
$
\end{ex}

\begin{prop}
\label{semiring-S0}Let $(L,\vee )$ be a join-semilattice and fix $c\in L.$

\begin{enumerate}
\item We have a structure of a \emph{non-unital} (additively idempotent)
demiring on $D_{\vee }^{c}(L):=(L,\vee ,\cdot _{c})$ under the following
operations%
\begin{equation*}
x+y:=x\vee y,\qquad x\cdot _{c}y=c\quad \text{for all }x,y\in L.
\end{equation*}

\item $D_{\vee }^{c}(L):=(L;\vee ,c;\cdot _{c})$ is a \emph{non-unital}
hemiring if and only if $L$ is bounded with bottom $c=0_{L}.$

\item Let $L$ be bounded with top element $1\neq 0$ and define%
\begin{equation}
x+y:=x\vee _{L}y,\quad \text{ }x\cdot _{\vee }y:=%
\begin{cases}
x, & y=1, \\ 
y, & x=1, \\ 
0, & x\neq 1\neq y.%
\end{cases}
\label{plus-dot-S0}
\end{equation}%
Then $S_{\vee }(L):=(L;\vee ,0;\cdot _{\vee },1)$ is a semiring if and only
if $1$ is join-irreducible.
\end{enumerate}
\end{prop}

\begin{Beweis}
The addition is associative, commutative, and idempotent because it is the
join of a join-semilattice. The multiplication is clearly commutative.

\begin{enumerate}
\item The multiplication is constant, whence associative. For every $%
x,y,z\in L$, we have%
\begin{equation*}
x\cdot _{c}(y+z)=c=c\vee c=(x\cdot _{c}y)+(x\cdot _{c}z).
\end{equation*}%
If $e$ were a multiplicative identity, then $x=e\cdot _{c}x=c$ for every $%
x\in L$, contrary to the assumption $\left\vert L\right\vert \geq 2.$

\item Obvious.

\item Let $L$ be bounded with top element $1\neq 0.$ We verify
associativity. Let $x,y,z\in L$. If $x=1$, then 
\begin{equation*}
(x\cdot _{\vee }y)\cdot _{\vee }z=y\cdot _{\vee }z=1\cdot _{\vee }(y\cdot
_{\vee }z)=x\cdot _{\vee }(y\cdot _{\vee }z).
\end{equation*}%
If $y=1$, then%
\begin{equation*}
(x\cdot _{\vee }y)\cdot _{\vee }z=x\cdot _{\vee }z=x\cdot _{\vee }(y\cdot
_{\vee }z).
\end{equation*}%
If $z=1$, then%
\begin{equation*}
(x\cdot _{\vee }y)\cdot _{\vee }z=x\cdot _{\vee }y=x\cdot _{\vee }(y\cdot
_{\vee }z).
\end{equation*}%
If none of $x,y,z$ is $1$, then $xy=yz=0$, and%
\begin{equation*}
(x\cdot _{\vee }y)\cdot _{\vee }z=0\cdot _{\vee }z=0=x\cdot _{\vee }0=x\cdot
_{\vee }(y\cdot _{\vee }z),
\end{equation*}%
where we used that $0$ is absorbing for multiplication. Hence multiplication
is associative. Clearly $1$ is its identity and $0$ is absorbing. So, it
remains to prove that the distributivity holds if and only if $1$ is
join-irreducible.

$(\Longrightarrow )$ Assume that $S_{\vee }(L)$ is a semiring. Suppose that $%
1$ is \emph{not} join-irreducible. Then there exist $y,z<1$ such that $y\vee
_{L}z=1$. Neither $y$ nor $z$ is $0$, since $0\vee _{L}z=z$ and $0\vee
_{L}y=y$. In particular, $y\neq 0$. Since $y,z\neq 1$, the multiplication
rule gives%
\begin{equation*}
y\cdot _{\vee }(y+z)=y\cdot _{\vee }1=y,\qquad y\cdot _{\vee }y+y\cdot
_{\vee }z=0\vee _{L}0=0.
\end{equation*}%
As $y\neq 0$, distributivity fails, a contradiction.

$(\Longleftarrow )$ Assume that the top element $1$ is join-irreducible.

\textbf{Claim: }By the commutativity of the multiplication $\cdot _{\vee },$
we need only to prove%
\begin{equation}
x\cdot _{\vee }(y+z)=(x\cdot _{\vee }y)+(x\cdot _{\vee }z)\text{ for all }%
x,y,z\in L.  \label{dist-over-L}
\end{equation}

\textbf{Case I: }$x=1.$ In this case, the identity (\ref{dist-over-L}) is
immediate.

\textbf{Case II: }$x\neq 1$ and $y=1$ (or $z=1$). Assume that $y=1.$ In this
case, we have $x\cdot _{\vee }(y+z)=x\cdot _{\vee }(1\vee z)=x\cdot _{\vee
}1=x.$

\begin{itemize}
\item If $z=1,$ then%
\begin{equation*}
(x\cdot _{\vee }y)+(x\cdot _{\vee }z)=x\cdot _{\vee }1+x\cdot _{\vee
}1=x\vee x=x=x\cdot _{\vee }1=x.
\end{equation*}

\item If $z\neq 1,$ then%
\begin{equation*}
(x\cdot _{\vee }y)+(x\cdot _{\vee }z)=(x\cdot _{\vee }1)\vee (x\cdot _{\vee
}z)=x\vee 0=x.
\end{equation*}
\end{itemize}

In either subcases, the distributivity condition (\ref{dist-over-L})\ holds.

\textbf{Case III: }$x,y,z\neq 1$. By the join-irreducibility of $1$, we have 
$y\vee z\neq 1$ and so%
\begin{equation*}
x\cdot _{\vee }(y+z)=x\cdot _{\vee }(y\vee z)=0=0\vee 0=(x\cdot _{\vee
}y)+(x\cdot _{\vee }z).
\end{equation*}

In all cases, the distributivity holds and therefore $S_{\vee }(L)$ is an
additively idempotent semiring.$\blacksquare $
\end{enumerate}
\end{Beweis}

\begin{rem}
\label{non-iso-hem}Let $(L,\vee ,0)$ be a join-semilattice with a
join-irreducible $1\neq 0.$ The hemirings $(S_{\vee }(L);\vee ,0;\cdot
_{\vee })$ and $(H_{\vee }(L);\vee ,0;\cdot _{0})$ are \emph{%
non-isomorphic.\ }In fact, there is even \emph{no} surjective morphisms of
hemirings from $S_{\vee }(L)$ to $H_{\vee }(L):$ Suppose $\varphi :S_{\vee
}(L)\longrightarrow H_{\vee }(L)$ were a surjective morphism of hemirings.
Notice that%
\begin{equation*}
\varphi (1)=\varphi (1\cdot _{S_{\vee }(L)}1)=\varphi (1)\cdot _{H_{\vee
}(L)}\varphi (1)=0.
\end{equation*}%
Let $h\in H_{\vee }(L)$ and pick $s\in \varphi ^{-1}(h)\subseteq S_{\vee
}(L).$ Then%
\begin{equation*}
\begin{tabular}{lllll}
$h$ & $=$ & $\varphi (s)$ & $=$ & $\varphi (s\cdot _{S_{\vee }(L)}1)$ \\ 
& $=$ & $\varphi (s)\cdot _{H_{\vee }(L)}\varphi (1)$ & $=$ & $\varphi
(s)\cdot _{H_{\vee }(L)}0$ \\ 
& $=$ & $0.$ & $=$ & 
\end{tabular}%
\end{equation*}%
So, $H_{\vee }(L)=\{0\},$ a contradiction to $\left\vert L\right\vert \geq
2.\blacksquare $
\end{rem}

\begin{ex}
\label{C2-join-S}Consider the two-element chain $C_{2}:\bot <\top .$

\begin{enumerate}
\item Fixing $c=1,$ the operations of $D_{\vee }^{1}(C_{2})$ are given by%
\begin{equation*}
\begin{tabular}{ccccc}
\begin{tabular}{ccc}
$\vee $ & $0$ & $1$ \\ 
$0$ & $0$ & $1$ \\ 
$1$ & $1$ & $1$%
\end{tabular}
&  &  &  & 
\begin{tabular}{ccc}
$\cdot _{1}$ & $0$ & $1$ \\ 
$0$ & $1$ & $1$ \\ 
$1$ & $1$ & $1$%
\end{tabular}%
\end{tabular}%
\end{equation*}%
Clearly, $D_{\vee }^{1}(C_{2})$ is a non-unital demiring.

\item Choosing $c=0,$ the operations of $H_{\vee }(C_{2}):=D_{\vee
}^{0}(C_{2})$ are%
\begin{equation*}
\begin{tabular}{ccccc}
\begin{tabular}{ccc}
$\vee $ & $0$ & $1$ \\ 
$0$ & $0$ & $1$ \\ 
$1$ & $1$ & $1$%
\end{tabular}
&  &  &  & 
\begin{tabular}{ccc}
$\cdot _{0}$ & $0$ & $1$ \\ 
$0$ & $0$ & $0$ \\ 
$1$ & $0$ & $0$%
\end{tabular}%
\end{tabular}%
\end{equation*}%
Then $H_{\vee }(C_{2})\cong \mathfrak{B}$ (the two-element \emph{Boolean
hemiring}) which is clearly non-unital.

\item By the definition of $S_{\vee }(C_{2}),$ it has the following
operations%
\begin{equation*}
\begin{tabular}{ccccc}
\begin{tabular}{ccc}
$\vee $ & $0$ & $1$ \\ 
$0$ & $0$ & $1$ \\ 
$1$ & $1$ & $1$%
\end{tabular}
&  &  &  & 
\begin{tabular}{ccc}
$\cdot _{\vee }$ & $0$ & $1$ \\ 
$0$ & $0$ & $0$ \\ 
$1$ & $0$ & $1$%
\end{tabular}%
\end{tabular}%
\end{equation*}%
and so $S_{\vee }(C_{2})\cong \mathbb{B}$ (the \textbf{two-element Boolean
semiring}).$\blacksquare $
\end{enumerate}
\end{ex}

\begin{rem}
\label{bounded}Let $(T,\vee )$ be a join-semilattice and assume that $T$ has 
\emph{no top} or its top $1$ is \emph{not} join-irreducible. We construct a
bounded join-semilattice $L$ whose \emph{new} greatest element is \emph{%
join-irreducible}. Retain the order of $T$ and adjoin to it a \emph{new}
element $1^{\sharp }$ that is \emph{not in} $T.$ If $T$ has a least element $%
0,$ declare $x<1^{\sharp }$ for all $x\in T.$ If $T$ has no least element,
add a \emph{new} element $\mathbf{0},$ set $L:=T\cup \{\mathbf{0},$ $%
1^{\sharp }\}$ and declare%
\begin{equation*}
\mathbf{0}<x<1^{\sharp }\qquad \text{for every }x\in T\text{ and }\mathbf{0}%
<1^{\sharp }.
\end{equation*}%
Then $L$ is a bounded join-semilattice, with least element $\mathbf{0}$ and
a \emph{join-irreducible} greatest element $1^{\sharp }.$
\end{rem}

The following example illustrates that the \emph{join-irreducibility}
condition (\ref{top-join-prime}) is essential for $S_{\vee }(L)\ $to be a
semiring.

\begin{ex}
Consider the four-element Boolean lattice%
\begin{equation*}
B_{2}:e<a,\text{ }b<t,\text{ }a\parallel b.
\end{equation*}%
Notice that $t$ is the top of $B_{2}$ but $t$ is \emph{not}
join-irreducible: we have $a\vee b=t$ with $a,b<t.$ Let's apply the
multiplication formula (\ref{plus-dot-S0}) with $1=t$ and $0=e.$ Since $%
a+b=a\vee _{L}b=t$ and $a,b<t$, one obtains 
\begin{equation*}
a\cdot _{\vee }(a+b)=a\cdot _{\vee }t=a\neq e=e+e=a\cdot _{\vee }a+a\cdot
_{\vee }b.
\end{equation*}%
Thus \emph{distributivity fails} and this construction does \emph{not}
define a semiring on $B_{2}.$

Adding a genuinely \emph{new} greatest element $1^{\sharp }$ \emph{not} in $%
B_{2}$ avoids this obstruction. On the other hand, the pre-existing least
element $e$ can safely be reused by setting $\mathbf{0}=e,$ as in the
construction above.$\blacksquare $
\end{ex}

\begin{rem}
The idea of adjoining a new top element to a semilattice to construct a
semiring from a bounded join-semilattice should be well known. Although we
could not locate a formal publication exposing this construction in the way
we did, we discovered while revising the manuscript a public answer on \emph{%
MathOverflow} by Keith Kearnes \cite{MO20025} (posted 1 November 2025) which
gave the same basic construction of the semiring $S_{\vee }(L)$ associated
to a join-semilattice $L.$
\end{rem}

\begin{defn}
(cf. \cite[Ch. 20, Ch. 21]{Gol1999}) Let $(S,+,\cdot )\ $be a demiring which
is also a poset $(S,\leq )$ and fix $c\in S.$ We say that $S$ is a

\begin{enumerate}
\item $c$-\textbf{partially-ordered demiring} iff for all $%
s_{1},s_{2},s,t\in S$ with $t\geq c$ we have 
\begin{equation}
s_{1}\leq s_{2}\Longrightarrow s+s_{1}\leq s+s_{2},\quad ts_{1}\leq
ts_{2},\quad s_{1}t\leq s_{2}t.  \label{POH}
\end{equation}

\item \textbf{lattice-ordered-demiring} iff $S$ also has the structure of a
lattice $(S,\wedge ,\vee )$ such that for all $s,t\in S$ we have%
\begin{equation}
s+t=s\vee t,\qquad s\cdot t\leq s\wedge t.  \label{lat-ordered}
\end{equation}
\end{enumerate}
\end{defn}

\begin{ex}
Let $(L,\vee )$ be a join-semilattice and fix $c\in L.$ The induced demiring 
$D_{\vee }^{c}(L)$ is \emph{additively idempotent,} and so it has the
natural partial order $x\leq y$ iff $x+y=y$ (which is precisely the order of 
$L$). Since the multiplication is constant, $S_{\vee }^{c}(L)$ is clearly a $%
c$-partially-ordered demiring.

If $L$ is bounded with bottom $0,$ then $t\geq 0$ for all $t\in H_{\vee }(L)$
and $H_{\vee }(L)$ is a \emph{partially-ordered hemiring} in the sense of 
\cite[Ch. 20]{Gol1999}.$\blacksquare $
\end{ex}

\begin{ex}
Let $(L,\wedge ,\vee )$ be a lattice and fix $c\in L.$ Then $D_{\vee
}^{c}(L) $ is clearly a lattice-ordered demiring if and only if $L$ has a
bottom $0$ and $c=0$ (in this case, $H_{\vee }(L)=D_{\vee }^{0}(L)$ is a
lattice-ordered-hemiring). If $(L,\wedge ,\vee ,0,1)$ is bounded with
join-irreducible $1\neq 0,$ then each product $x\cdot y$ in $S_{\vee }(L)$
is either $0$ or the meet (when $1\in \{x,y\}$), whence $s\cdot _{\vee
}t\leq s\wedge t$ for all $s,t$ and $S_{\vee }(L)$ is a \emph{%
lattice-ordered semiring}.$\blacksquare $
\end{ex}

\begin{rem}
Our semiring induced by a join-semilattice\emph{\ }$S_{\vee }(L)$ is not
necessarily a \emph{semilattice-ordered semiring} (abbreviated \textbf{%
sl-semiring}) in the sense of Mala and Rao \cite{RM2008} (cf. \cite{MR2016}, 
\cite{Rao1981}). Besides a semilattice order, an sl-semiring has an extra
subtraction/residuation operation satisfying additional axioms. The semiring 
$S_{\vee }(L)$ does not necessarily have this extra operation.
\end{rem}

\begin{ex}
\label{not-sl-sem}Consider the \emph{diamond}%
\begin{equation*}
M_{3}:0<a,b,c<t.
\end{equation*}%
Adjoin a fresh top element $1^{\sharp }$ (with $t<1^{\sharp }$) and set $%
L:=M_{3}\cup \{1^{\sharp }\}.$ The new top $1^{\sharp }$ is clearly
join-irreducible and so $S_{\vee }(L)$ is a semiring by Proposition \ref%
{semiring-S0}.

\textbf{Claim: }$S_{\vee }(L)$ is \emph{not} an sl-semiring.\ 

A necessary condition of an sl-semiring is that for every $a,b$ there be a
least $x$ satisfying $a\leq b\vee x$ (cf. \cite[Definition 2.4]{MR2016}).

For target $a$ and fixed $b,$ both $x=a$ and $x=c$ satisfy $a\leq b\vee x$
(since $b\vee a=t=b\vee c$). A least such $x$ would have to satisfy $x\leq a$
and $x\leq c,$ hence $x\leq a\wedge c=0.$ But $b\vee 0=b$ is not above $a.$
Therefore no least residual exists, and $S_{\vee }(L)$ fails this necessary
sl-semiring axiom.$\blacksquare $
\end{ex}

\begin{prop}
\label{ideals-join}Let $(L,\vee )\ $be a join-semilattice, fix $c\in L$ and
consider the lattice of order ideals of $U_{c}:$%
\begin{equation}
\func{Idl}(U_{c})=\{D\subseteq \lbrack c\uparrow )\mid D\text{ is down-set
closed under}\vee \}.  \label{Idl(Uc)}
\end{equation}

\begin{enumerate}
\item The collection of \emph{all} ideals of $D_{\vee }^{c}(L)$ is%
\begin{equation}
\mathfrak{I}(D_{\vee }^{c}(L))=\{I\subseteq L\text{ }|\text{ }c\in I\text{
and }I\text{ closed under }\vee \}.  \label{I(S0(L))}
\end{equation}

\item The subtractive closure of an ideal $I$ of $D_{\vee }^{c}(L)$ is%
\begin{equation*}
\overline{I}=\{x\in L\text{ }|\text{ }x\leq h\text{ for some }h\in
I\}=\bigcup_{h\in I}(\downarrow h].
\end{equation*}

\item An ideal $I$ of $D_{\vee }^{c}(L)$ is subtractive if and only if $I$
is an order ideal of $L$ containing $c.$

\item There is a canonical isomorphism of complete lattices%
\begin{equation}
\mathfrak{K}(D_{\vee }^{c}(L))\cong \func{Idl}(U_{c}).  \label{can-iso-join}
\end{equation}

\item $\mathfrak{K}(D_{\vee }^{c}(L))$ is a complete algebraic lattice with
bottom $(\downarrow c]$ and top $L,$ whose compact elements are precisely
the principal order ideals $\{(\downarrow d]$ $|$ $c\leq d\}.$
\end{enumerate}
\end{prop}

\begin{Beweis}
First of all, notice that for any $c\in L,$ the upper cone $%
U_{c}:=[c\uparrow )$ is a bounded semilattice with bottom $c,$ whence $\func{%
Idl}(U_{c})$ is a complete algebraic lattice by Lemma \ref{join-sem-idl} (2).

\begin{enumerate}
\item Let $I$ be an ideal of $D_{\vee }^{c}(L).$ Then $I$ is closed under $%
\vee $ since the join is the addition of the demiring, and for any $a\in I$
we have $c=a\cdot a\in I.$

Conversely, let $I\subseteq L$ contain $c$ and be closed under $\vee .$ Then 
$I$ is an additive subsemigroup of $L.$ Moreover, if $s\in L$ and $x\in I$,
then $sx=c\in I.$ So, $I$ is an ideal of $D_{\vee }^{c}(L).$

\item Let $I$ be an ideal of $D_{\vee }^{c}(L).$ If $x\in \overline{I}$,
then $x\leq x\vee _{L}a:=x+a=\widetilde{a}$ for some $a,\widetilde{a}\in I$.
Conversely, if $x\leq h$ for some $h\in I$, then $x+h:=x\vee _{L}h=h$ and so 
$x\in \overline{I}.$

\item Let $I\in \mathfrak{K}(D_{\vee }^{c}(L)).$ If $x\in I$ and $y\leq x,$
then $y+x=y\vee _{L}x=x\in I$ and so, $y\in \overline{I}=I$ (since $I$ is
subtractive). Thus $I$ is downward closed. Moreover, $c\in I$ and closed
under joins because $I$ is an ideal, so $I$ is a non-empty order ideal of $L$
that contains $c.$

Conversely, let $A$ be an order ideal of $L$ containing $c.$ Then $I$ is
closed under joins and so an ideal of $D_{\vee }^{c}(L)$ by (1). If $%
x,x+y\in A$, then $y\leq x\vee _{L}y=x+y$ and the downward closure of $A$
gives $y\in A$, so $A$ is subtractive.

\item For any $A\subseteq L$ set%
\begin{equation*}
(\downarrow A]:=\{x\in L\mid x\leq a\text{ for some }a\in A\}.
\end{equation*}%
Consider the maps%
\begin{eqnarray*}
\Phi &:&\mathfrak{K}(D_{\vee }^{c}(L))\longrightarrow \func{Idl}(U_{c}),%
\text{ }I\longmapsto I\cap U_{c}; \\
\Psi &:&\func{Idl}(U_{c})\longrightarrow \mathfrak{K}(D_{\vee }^{c}(L)),%
\text{ }A\longmapsto (\downarrow A].
\end{eqnarray*}%
\textbf{Claim:}\ $\Phi $ and $\Psi $ are well define.

If $I$ is a subtractive ideal of $D_{\vee }^{c}(L),$ then $I\cap U_{c}$ is
nonempty as it contains $c,$ downward closed in $U_{c},$ and join-closed. On
the other hand, if $A\in \func{Idl}(U_{c}),$ then $c\in (\downarrow A]$
which is downward closed and join closed: if $x\leq a$ and $y\leq b$ for
some $a,b\in A,$ then $x\vee y\leq a\vee b\in A.$

For $I\in \mathfrak{K}(D_{\vee }^{c}(L))$ and $x\in I$, the element $x\vee c$
belongs to $I\cap U_{c}$ and $x\leq x\vee c.$ Hence $\Psi (\Phi
(I))=(\downarrow (I\cap U_{c})]=I.$

On the other hand, for $A\in \func{Idl}(U_{c})$, if $u\in (\downarrow A]\cap
U_{c},$ then $u\leq a$ for some $a\in A$ and the down-closure of $A$ inside $%
U_{c}$ gives $u\in A$. Therefore $(\downarrow _{L}A]\cap U_{c}=A.$

So, $\Phi $ and $\Psi $ are inverse order isomorphisms, hence complete
lattice isomorphisms.

\item Any ideal lattice of a bounded join-semilattice is complete algebraic
with compact elements the principal order ideals. Applying this to $U_{c}$
gives the final claims. Explicitly, meets are intersections and 
\begin{equation*}
A\vee B=\{u\in \lbrack c\uparrow )\text{ }\mid \text{ }u\leq a\vee b\text{
for some }a\in A,\ b\in B\}.
\end{equation*}%
The principal order ideal $(_{U_{c}}\downarrow u]$ is sent by $\Psi $ to $%
(_{L}\downarrow u].\blacksquare $
\end{enumerate}
\end{Beweis}

\begin{cor}
\label{ideals-join-S}Let $(L,\vee ,0)$ be a bounded join-semilattice with
join-irreducible $1\neq 0.$

\begin{enumerate}
\item The collection of \emph{all} ideals of $S_{\vee }(L)$ is%
\begin{equation}
\mathfrak{I}(S_{\vee }(L))=\{L\}\cup \{I\subseteq L\backslash \{1\}\text{ }|%
\text{ }0\in I,\text{ }I\text{ closed under }\vee \}.  \label{ideal-S-jn}
\end{equation}

\item There is a canonical isomorphism of complete lattices%
\begin{equation}
\mathfrak{K}(S_{\vee }(L))\cong \func{Idl}(L).  \label{k-ideal-sem-L}
\end{equation}
\end{enumerate}
\end{cor}

\begin{Beweis}
Let $L$ be a bounded join-semilattice with join-irreducible top $1\neq 0$.

\begin{enumerate}
\item This follows from Proposition \ref{ideals-join} by choosing $c=0.$

\item The subtractive closure depends only on the addition, so $\mathfrak{K}%
(S_{\vee }(L))\subseteq \mathfrak{K}(H_{\vee }(L))\overset{\text{Proposition %
\ref{ideals-join} (4)}}{=}\func{Idl}(L).$ Conversely, let $D\in \func{Idl}%
(L).$ If $1\in D,$ then $D=L.$ If $1\notin D,$ then $D\subseteq L\backslash
\{1\}$ and $0\in D.$ In both cases $D$ is an ideal of $S_{\vee }(L)$ by (1).
Moreover, $D$ is a subtractive ideal because it is down-closed. Hence, $%
\mathfrak{K}(S_{\vee }(L))=\mathfrak{K}(H_{\vee }(L))=\func{Idl}(L),$ and
the lattice operations agree since $\wedge _{\mathfrak{K}}=\cap $ and $I\vee
_{\mathfrak{K}}J=\overline{I+J}$ involve only the addition.$\blacksquare $
\end{enumerate}
\end{Beweis}

\begin{rem}
\label{ideal-join-1}Let $(L,\vee ,0)$ be a bounded join-semilattice with
join-irreducible $1\neq 0.$ The exclusion of $1$ in the second family of \ref%
{ideal-S-jn} is essential. For example, if $\left\vert L\right\vert \geq 3$
then $I=\{0,1\}$ contains $0$ and is join-closed but $I$ is \emph{not} an
ideal for any $a\in L\backslash I$ we have $a\cdot _{\vee }1=a\notin
\{0,1\}.\blacksquare $
\end{rem}

Using duality and applying Propositions \ref{semiring-S0} and \ref%
{ideals-join}, we get

\begin{prop}
\label{semiring-meet}Let $(M,\wedge )$ be a meet-semilattice.

\begin{enumerate}
\item Fix $c\in M.$ We have a structure of a \emph{non-unital} (additively
idempotent) demiring on $D_{\wedge }^{c}(M):=(M,\wedge ,\cdot _{c})$ under
the following operations%
\begin{equation*}
x+y:=x\wedge y,\qquad x\cdot _{c}y=c\quad \text{for all }x,y\in M.
\end{equation*}%
Moreover, $(M^{\mathrm{op}},\wedge ^{\mathrm{op}})$ is a join-semilattice and%
\begin{equation*}
D_{\wedge }^{c}(M)\cong D_{\wedge ^{\mathrm{op}}}^{c}(M^{\mathrm{op}}),\quad 
\mathfrak{K}(D_{\wedge }^{c}(M))\cong \func{Fil}(D_{c}).
\end{equation*}

\item $H_{\wedge }(M):=(M;\wedge ,c;\cdot _{c})$ is a \emph{non-unital}
hemiring if and only if $M$ is bounded with top $c=\top .$ In this case, $%
(M^{\mathrm{op}},\wedge ^{\mathrm{op}},\top )$ is a bounded join-lattice
with least element $0=\top $ and%
\begin{equation*}
H_{\wedge }(M)\cong H_{\wedge ^{\mathrm{op}}}(M^{\mathrm{op}}),\quad 
\mathfrak{K}(H_{\wedge }(M))\cong \func{Fil}(M).
\end{equation*}

\item Let $(M,\wedge ,\top )$ be bounded with bottom $\bot \neq \top $ and
define%
\begin{equation}
x+y:=x\wedge y,\quad \text{ }x\cdot _{\wedge }y:=%
\begin{cases}
x, & y=\bot , \\ 
y, & x=\bot , \\ 
\top , & x\neq \bot \neq y.%
\end{cases}%
\end{equation}%
Then $S_{\wedge }(M):=(M;\wedge ,\top ;\cdot _{\wedge },\bot )$ is a
semiring if and only if the bottom $\bot $ is meet-irreducible. In this
case, $(M^{\mathrm{op}},\wedge ^{\mathrm{op}},\top )$ is a bounded
join-lattice with a join-irreducible greatest element $1=\bot \neq 0=\top $
and 
\begin{equation*}
S_{\wedge }(M)\cong S_{\wedge ^{\mathrm{op}}}(M^{\mathrm{op}}),\quad 
\mathfrak{K}(S_{\wedge }(M))\cong \func{Fil}(M).
\end{equation*}
\end{enumerate}
\end{prop}

The following result follows directly from the definitions.

\begin{prop}
\label{coideal-join}Let $(L,\vee )\ $be a join-semilattice and fix $c\in L.$

\begin{enumerate}
\item $\mathcal{A}(D_{\vee }^{c}(L))=\mathfrak{U}(L).$

\item $\mathcal{C}(D_{\vee }^{c}(L))=\{C\in \mathfrak{U}(L)$ $\mid $ $c\in
C\}=\mathcal{B}(D_{\vee }^{c}(L)).$
\end{enumerate}
\end{prop}

\begin{prop}
\label{j-semiring-local}Let $(L,\vee ,0,1)$ be a bounded join-semilattice
with a join-irreducible top $1\neq 0.$ The semiring $S_{\vee }(L)$ is local
and, moreover, its unique maximal ideal $\mathfrak{m}_{\vee }:=L\backslash
\{1\}$ is subtractive.
\end{prop}

\begin{Beweis}
Consider ideal $\mathfrak{m}_{\vee }:=L\backslash \{1\}.$ Then $\mathfrak{m}%
_{\vee }$ is closed under addition (joins) since $1$ is join-irreducible.
Moreover, $\mathfrak{m}_{\vee }$ is clearly stable under multiplication.
Clearly $\mathfrak{m}_{\vee }$ is maximal and contains all proper ideals of $%
S_{\vee }(L).$ Notice that $\mathfrak{m}_{\vee }$ is subtractive: if $x\in 
\mathfrak{m}_{\vee }$ and $x\vee y\in \mathfrak{m}_{\vee }$, then $y\neq 1$,
so $y\in \mathfrak{m}_{\vee }.\blacksquare $
\end{Beweis}

The following example shows that a demiring/hemiring induced by a
semilattice might have\emph{\ no} maximal ideals.

\begin{ex}
\label{dyadic}Let $L$ be the Boolean algebra of \emph{finite} unions of all
dyadic half-open subintervals%
\begin{equation*}
\lbrack \frac{k}{2^{n}},\frac{k+1}{2^{n}}),\text{ }\qquad n\geq 0,\qquad
0\leq k<2^{n}
\end{equation*}%
of $[0,1)$ with join $A\vee B:=A\cup B$ and bottom element $\varnothing .$%
Every \emph{non-empty} $A\in L$ can be written as $A=A_{1}\cup A_{2}$ with $%
A_{1},A_{2}\varsubsetneqq A$ (bisect one dyadic interval in a representation
of $A$).

\textbf{Claim: }$H_{\vee }(L)$ has \emph{no} maximal ideals. Suppose that $%
\mathfrak{m}$ were a maximal ideal. Choose $A\notin \mathfrak{m}$. Since
multiplication is zero, ideals are exactly the union-closed subsets
containing $\varnothing .$ The ideal $I$ generated by $\mathfrak{m}\cup A$
consists of elements of the form $B$ or $B\cup A$ with $B\in \mathfrak{m}.$
Since $m\varsubsetneqq I\subseteq H_{\vee }(L),$ the maximality of $%
\mathfrak{m}$ forces $I=L.$ Choose $A_{1},A_{2}\in L$ with $A=A_{1}\cup
A_{2} $ and $A_{1},A_{2}\varsubsetneqq A.$ Neither $A_{1}$ nor $A_{2}$ could
be of the form $A\cup B.$ Since $I=L,$ we have $A_{1},A_{2}\in \mathfrak{m}$
and so $A=A_{1}\cup A_{2}\in \mathfrak{m}$, a contradiction. So, $H_{\vee
}(L)$ has no maximal ideals.$\blacksquare $
\end{ex}

The following example shows that a demiring/hemiring induced by a
semilattice might have a maximal ideal which is \emph{not} prime.

\begin{ex}
\label{max-not-prime}Consider the chain $C_{2}:0<a.$ The hemiring $H_{\vee
}(C_{2})$ has only one proper ideal $\mathfrak{m}:=\{0\}$ which is clearly
maximal but \emph{not} prime (as $a^{2}=0\in \mathfrak{m}$ but $a\notin 
\mathfrak{m}$).
\end{ex}

\begin{lem}
\label{downset-hom}Let $(L,\vee )$ be a join-semilattice and consider%
\begin{equation}
f:L\longrightarrow \func{Idl}(L),\qquad a\longmapsto (\downarrow a].
\label{f:L-Idl(L)}
\end{equation}

\begin{enumerate}
\item $f$ is an injective morphism of join-semilattices and an order
embedding. Moreover, if $L$ is bounded with bottom $0$ (has top $1$), then $%
\func{Idl}(L)$ is bounded with bottom $\{0\}$ (has top $f(1)=L$).

\item If $(L,\wedge ,\vee ,0,1)$ is a bounded lattice, then $f$ is an
injective morphisms of bounded lattices. In particular, $L$ is isomorphic to
the bounded sublattice $\mathfrak{D}_{p}(L)\subseteq \func{Idl}(L)$ of
principal ideals of $L.$
\end{enumerate}
\end{lem}

\begin{Beweis}
We first observe that $f$ is well defined since $f(a)=(\downarrow a]$ is an
order ideal for every $a\in L.$ For all $a,b\in L,$ we have $a\leq
b\Longleftrightarrow (\downarrow a]\subseteq (\downarrow b]$ whence $f$ is
an injective order embedding. Moreover, if $L$ is bounded with bottom $0,$
then $f(0)=(\downarrow 0]=\{0\}$ is clearly bottom in $\func{Idl}(L).$ If $L$
has a top $1,$ then $f(1)=(\downarrow 1]=L$ which is the top of $\func{Idl}%
(L).$

\begin{enumerate}
\item \textbf{Claim:}\ $f$ preserves joins. Recall that the join of order
ideals given in \ref{j-idl}. In particular, for any $a,b\in L,$ we have%
\begin{eqnarray*}
f(a)\vee _{\func{Idl}(L)}f(b) &=&\{x\in L\text{ }|\text{ }x\leq u\vee v\text{
for some }u\leq a,\ v\leq b\} \\
&\subseteq &\{x\in L\text{ }|\text{ }x\leq a\vee b\}=f(a\vee b)
\end{eqnarray*}%
Conversely, take any $x\in f(a\vee b)=(\downarrow a\vee b]$. Then $x\leq
a\vee b.$ Since $a\in f(a)$ and $b\in f(b)$, it follows by (\ref{j-idl})
that $x\in f(a)\vee _{\func{Idl}(L)}f(b).$ Therefore, we have $f(a\vee
_{L}b)\subseteq f(a)\vee _{\func{Idl}(L)}f(b).$ Combining the two
inclusions, we obtain%
\begin{equation*}
f(a\vee _{L}b)=f(a)\vee _{\func{Idl}(L)}f(b).
\end{equation*}

\item \textbf{Claim:}\ $f$ preserves meets. Let $a,b\in L$. Since the meet
in $\func{Idl}(L)$ is intersection, we have%
\begin{equation*}
f(a)\wedge _{\func{Idl}(L)}f(b)=f(a)\cap f(b)=(\downarrow a]\cap (\downarrow
b]=(\downarrow a\wedge _{L}b]=f(a\wedge _{L}b).
\end{equation*}%
Since $f$ preserves joins by (1) and meets by the argument above, we
conclude that $f$ is a morphism of bounded lattices. Since $f$ is injective,
we have $L\cong \mathfrak{D}_{p}(L)$ which is a bounded sublattice of $\func{%
Idl}(L).\blacksquare $
\end{enumerate}
\end{Beweis}

\begin{ex}
Consider the set of whole numbers $\mathbb{W}=\{0,1,2,\cdots \}$ with its
usual order and join $x\vee y=\max \{x,y\}.$ Then $(\mathbb{W},\max ;0)$ is
a bounded join-semilattice with bottom $0$ and \emph{no} top. Its order
ideals are exactly 
\begin{equation*}
\func{Idl}(\mathbb{W})=\{\mathbb{W}\}\cup \{0,1,\ldots ,n\}\text{ }|\text{ }%
n\in \mathbb{W}\}.
\end{equation*}%
The principal-ideal embedding is 
\begin{equation*}
n\longmapsto (\downarrow n]=\{0,1,\cdots ,n\}
\end{equation*}%
and it clearly preserves $0$ and binary joins, since $(\downarrow \max
\{m,n\}]=(\downarrow m]\vee (\downarrow n].$ The greatest element of $\func{%
Idl}(\mathbb{W})$ is $\mathbb{W}$, which is \emph{not} principal;
correspondingly, $\mathbb{W}$ itself has no greatest element. Thus Lemma \ref%
{downset-hom} gives an embedding preserving the existing bottom, and makes
no claim of preserving a nonexistent top.$\blacksquare $
\end{ex}

\begin{thm}
\label{mod-S-join}Let $(L,\vee )$ be a join-semilattice and fix $c\in L.$
The following are equivalent:

\begin{enumerate}
\item $\mathfrak{K}(D_{\vee }^{c}(L))$ is a modular lattice;

\item $\func{Idl}(U_{c})$ is a modular lattice;

\item $U_{c}$ is a modular join-semilattice.
\end{enumerate}
\end{thm}

\begin{Beweis}
The equivalence $(1\Longleftrightarrow 2)$ follows from the lattice
isomorphism (\ref{can-iso-join}). The equivalence $(2\Longleftrightarrow 3)\ 
$holds for any join-semilattice as stated in \cite[pp. 282--283]{Ruz2017} (a
proof for the dual case of meet-semilattices can be found in \cite{Cor1974}).%
$\blacksquare $
\end{Beweis}

\begin{thm}
\label{dist-S-join}Let $(L,\vee )$ be a join-semilattice and fix $c\in L.$
The following are equivalent:

\begin{enumerate}
\item $\mathfrak{K}(D_{\vee }^{c}(L))$ is a distributive lattice;

\item $\func{Idl}(U_{c})$ is a distributive lattice;

\item $U_{c}$ is a distributive join-semilattice.
\end{enumerate}
\end{thm}

\begin{Beweis}
The equivalence $(1\Longleftrightarrow 2)$ follows from the lattice
isomorphism (\ref{can-iso-join}). The equivalence $(2\Longleftrightarrow 3)\ 
$is a classical fact about join-semilattices (e.g. \cite[Lemma II, 5.1]%
{Gra1998}).$\blacksquare $
\end{Beweis}

\begin{thm}
\label{1-hash}Let $(L,\vee ,0)$ be a bounded join-semilattice.

\begin{enumerate}
\item If $L$ has a join-irreducible top $1\neq 0,$ then

\begin{enumerate}
\item $\mathfrak{K}(S_{\vee }(L))$ is a modular lattice $\Longleftrightarrow 
$ $L$ is a modular join-semilattice.

\item $\mathfrak{K}(S_{\vee }(L))$ is a distributive lattice $%
\Longleftrightarrow $ $L$ is a modular join-semilattice.
\end{enumerate}

\item If $L$ has no join-irreducible top and we adjoin $1^{\sharp }\notin L,$
then setting $L^{\sharp }:=L\cup \{1^{\sharp }\}$ we have

\begin{enumerate}
\item $\mathfrak{K}(S_{\vee }(L^{\sharp }))$ is a modular lattice $%
\Longleftrightarrow $ $L^{\sharp }$ is a modular join-semilattice $%
\Longleftrightarrow $ $L$ is a modular join-semilattice.

\item $\mathfrak{K}(S_{\vee }(L^{\sharp }))$ is a distributive lattice $%
\Longleftrightarrow $ $L^{\sharp }$ is a distributive join-semilattice $%
\Longleftrightarrow $ $L$ is a distributive join-semilattice.
\end{enumerate}
\end{enumerate}
\end{thm}

\begin{Beweis}
Let $(L,\vee ,0)$ be a bounded join-semilattice.

\begin{enumerate}
\item The result follows from the canonical isomorphism of lattices (\ref%
{k-ideal-sem-L}).

\item It is obvious that $L$ is a modular (distributive) join-semilattice\
if and only if $L^{\sharp }$ is a modular (distributive) join-semilattice.
The result follows from Theorem \ref{mod-S-join}, Theorem \ref{dist-S-join}
and the facts that%
\begin{equation*}
\mathfrak{K}(S_{\vee }(L))\cong \func{Idl}(L^{\sharp })=\func{Idl}(L)\cup
\{L^{\sharp }\}.\blacksquare
\end{equation*}
\end{enumerate}
\end{Beweis}

\begin{cor}
\label{hereditary}Let $(L,\vee )$ be a join-semilattice and fix $c\in L.$ If 
$L$ is modular (distributive), then $\mathfrak{K}(D_{\vee }^{c}(L))$ is
modular (distributive).
\end{cor}

\begin{Beweis}
By Proposition \ref{Uc-hereditary}, the upper cone $U_{c}$ inherits the
modularity (distributivity)\ from $L.$ The result follows now by applying
Theorem \ref{mod-S-join} (Theorem \ref{dist-S-join}) to the join-semilattice 
$U_{c}.\blacksquare $
\end{Beweis}

The converse in Corollary~\ref{hereditary} is false:

\begin{ex}
Consider the smallest non-modular lattice $N_{5}$ and let $t$ be its top.
Then $U_{t}=\{t\}$ and so $\mathfrak{K}(D_{\vee }^{t}(N_{5}))=\{N_{5}\}\cong 
\func{Idl}_{t}(N_{5})$ is the one-element lattice, which is trivially
modular (distributive).$\blacksquare $
\end{ex}

\subsection*{The lattices $\mathfrak{K}(S)$ is non-modular}

Theorem \ref{mod-S-join} provide families of counterexamples to the \textbf{%
wrong claim} that the lattice of subtractive ideals of any commutative
semiring is modular. Begin with any \emph{non-modular} bounded
join-semilattice $(L,\vee _{L},0_{L})$ with a join-irreducible top $%
1_{L}\neq 0_{L}$ (non-modular bounded meet-semilattice $(M,\wedge _{M},\top
_{M})$ with a meet-irreducible bottom $\bot _{M}\neq \top _{M}$) and
construct the associated additively idempotent\ semiring $S_{\vee }(L)$ ($%
S_{\wedge }(M)$). Then the lattice $\mathfrak{K}(S_{\vee }(L))$ ($\mathfrak{K%
}(S_{\wedge }(M))$)\ of subtractive ideals is \emph{non-modular}.

\bigskip

The so called \emph{subtractive ideals} (called also $k$\emph{-ideals}) of a
semiring were introduced by Henriksen \cite{Hen1958} and proved to be very
useful in investigating proper semirings. However, Henriksen announced (with
no proof) that the lattice $\mathfrak{K}(S)$ of all subtractive ideals of
any commutative semiring is modular. This \textbf{wrong claim} was repeated
in \cite[Proposition 2.6]{GD2024} citing \cite[Proposition 6.38]{Gol1999}
and \cite[Section 3]{SGBS} \cite[Theorem 2.4]{SA1992}. It is worth
mentioning that the cited results are correct but were proved in special
contexts: \cite[Proposition 6.38]{Gol1999} was proved for \emph{subtractive
semirings }and \cite[Theorem 2.4]{SA1992} which was proved for the
sublattice of \emph{full} subtractive ideals over \emph{additively inverse
semirings}.

\bigskip

While Theorem \ref{mod-S-join} gives us a machine to produce counterexamples
to this wrong claim, we provide next \emph{concrete examples} of
semilattice-induced semirings $S_{\vee }(L)$ for which $\mathfrak{K}(S_{\vee
}(L))$ is non-modular. We illustrate the non-modularity of its lattice of
subtractive ideals directly without referring to Theorem \ref{mod-S-join}.

\bigskip

The following example shows that $\mathfrak{K}(S_{\vee }(L))$ is \emph{%
non-modular} even when $L$ is a lattice (if $L$ is non-modular).

\begin{ex}
\label{N5+1}Let $N_{5}=\{e,a,b,c,t\}$ be the pentagon lattice with 
\begin{equation*}
e<a<c<t,\qquad e<b<t,\text{ }a\parallel b\text{ and }c\parallel b.
\end{equation*}%
Its meet and join satisfy, in particular,%
\begin{equation*}
a\wedge b=e,\quad c\wedge b=e,\quad a\vee b=t,\quad b\vee c=t.
\end{equation*}%
Set $L:=N_{5}\sqcup \{1^{\sharp }\},$ where $1^{\sharp }$ is a \emph{new}
greatest element strictly above every element of $N_{5}$. The least element
of $L$ is $e$ (the existing least element). Thus $(L,\wedge _{L},\vee
_{L},e,1^{\sharp })$ is a bounded lattice with bottom $e$ and top $1^{\sharp
}$. Moreover, the top $1^{\sharp }$ is join-irreducible: if $x,y<1^{\sharp }$%
, then $x,y\in N_{5}$, and $x\vee _{L}y\leq t<1^{\sharp }$. Hence $x\vee
_{L}y=1^{\sharp }$ forces $x=1^{\sharp }$ or $y=1^{\sharp }$. By Proposition~%
\ref{semiring-S0}, the algebra $S_{\vee }(L)$ is a commutative additively
idempotent semiring.

Apart from the argument above, it can be double checked directly using the
following addition and multiplication tables%
\begin{equation*}
\begin{array}{c|cccccc}
+ & e & a & b & c & t & 1^{\sharp } \\ \hline
e & e & a & b & c & t & 1^{\sharp } \\ 
a & a & a & t & c & t & 1^{\sharp } \\ 
b & b & t & b & t & t & 1^{\sharp } \\ 
c & c & c & t & c & t & 1^{\sharp } \\ 
t & t & t & t & t & t & 1^{\sharp } \\ 
1^{\sharp } & 1^{\sharp } & 1^{\sharp } & 1^{\sharp } & 1^{\sharp } & 
1^{\sharp } & 1^{\sharp }%
\end{array}%
\qquad 
\begin{array}{c|cccccc}
\cdot _{\vee } & e & a & b & c & t & 1^{\sharp } \\ \hline
e & e & e & e & e & e & e \\ 
a & e & e & e & e & e & a \\ 
b & e & e & e & e & e & b \\ 
c & e & e & e & e & e & c \\ 
t & e & e & e & e & e & t \\ 
1^{\sharp } & e & a & b & c & t & 1^{\sharp }%
\end{array}%
\end{equation*}%
that $S_{\vee }(L):=(S_{\vee }(N_{5}\sqcup \{1^{\sharp }\});\vee ,e;\cdot
_{\vee },1^{\sharp })$ is a commutative semiring.

\textbf{Claim: }$\mathfrak{K}(S_{\vee }(L))=\{O,$ $A,$ $B,$ $C,$ $N_{5},$ $%
L\},$ where%
\begin{equation*}
O:=\{e\},\qquad A:=\{e,a\},\qquad B:=\{e,b\},\qquad C:=\{e,a,c\}.
\end{equation*}%
Before we prove this directly, we note that this is consistent with
Proposition~\ref{ideals-join} (3) (which states that the subtractive ideals
are precisely the order ideals of $L$).

First of all, we notice that each listed set is indeed a subtractive ideal.
The sets $O,A,B,C$ are downward closed and join-closed by direct inspection.
Moreover, $N_{5}=(\downarrow t]$ the principal ideal generated by $t$ and $%
L=(\downarrow 1^{\sharp }]$ is the whole lattice; both are also downward
closed and join-closed.

Secondly, we verify that the list is exhaustive. Let $I$ be a subtractive
ideal of $S_{\vee }(L).$ If $1^{\sharp }\in I$, then $I=S_{\vee }(L)$ since $%
x=x\cdot _{\vee }1^{\sharp }\in I$ for every $x\in S_{\vee }(L).$ If $I$ is
a proper subtractive ideal and $t\in I$, then $I=N_{5}:\ $if $y\in N_{5},$
then $y+t=y\vee t=t\in I,$ whence $y\in \overline{I}=I.$ Suppose next that $%
t\notin I.$ Then $I$ cannot contain both $a$ and $b:$ Otherwise, $t=a\vee
_{L}b\in I.$ If $b,t\notin I$ and $c\in I,$ then $a+c=c\in I$, whence $a\in 
\overline{I}=I$ whence $I=C.$ If $c\notin I$, the only possibilities are $O$%
, $A$, and $B$. Thus the displayed list is exhaustive.

The inclusion order is 
\begin{equation*}
O<A<C<N_{5}<L,\text{ }O<B<N_{5}<L,\text{ }B\parallel A,\text{ }B\parallel C.
\end{equation*}%
In particular, we have%
\begin{equation*}
A\cap B=O,\qquad B\cap C=O,\qquad A\cap C=A.
\end{equation*}%
We now compute the \emph{adjusted join} of $A$ and $B$. Their ordinary ideal
sum is 
\begin{equation*}
A+B=\{x+y\text{ }|\text{ }x\in A,\ y\in B\}=\{e,a,b,t\},
\end{equation*}%
So, $A+B$ is not subtractive: $c+t=t\in A+B$, whence $c\in \overline{A+B}$
but $c\notin A+B.$ A direct computation shows that $\overline{A+B}=N_{5}$
(which is consistent with applying Proposition~\ref{ideals-join} (2)).

Notice that $A\subseteq C,$ but%
\begin{eqnarray*}
(A\vee _{\mathfrak{K}}B)\wedge _{\mathfrak{K}}C &=&\overline{A+B}\cap
C=N_{5}\cap C=C; \\
A\vee _{\mathfrak{K}}(B\wedge _{\mathfrak{K}}C) &=&A\vee _{\mathfrak{K}%
}O=A\neq C.
\end{eqnarray*}%
So, this is a direct failure of the modular identity. Thus $\mathfrak{K}%
(S_{\vee }(L))$ is \emph{non-modular}. This agrees with Theorem~\ref%
{mod-S-join} since $L$ is non-modular (as $L$ contains $N_{5}$ as a
sublattice, $L$ is non-modular by Lemma \ref{equiv-lat-semi}).$\blacksquare $
\end{ex}

\begin{rem}
\label{6-minimal}Let $(L,\vee ,0)$ be a bounded join-semilattice with
join-irreducible top element $1\neq 0.$ If $\left\vert L\right\vert \leq 5,$
then $T:=L\backslash \{1\}$ is a bounded join-semilattice and $\func{Idl}(L)=%
\func{Idl}(T)\cup \{L\}$ as an ordered set. Notice that $\left\vert
T\right\vert \leq 4.$

\textbf{Case I:\ }$T=B_{2}$ (the \emph{four-element Boolean lattice}), i.e.%
\begin{equation*}
B_{2}:e\leq a\leq t,\qquad e\leq b\leq t,\qquad a\parallel b.
\end{equation*}%
In this case, we have%
\begin{equation*}
\func{Idl}(B_{2}):\{e\}<\{e,a\},\{e,b\}<B_{2},
\end{equation*}%
which is isomorphic to the Boolean lattice $B_{2}.$

\textbf{Case II:\ }$T$ is a chain, whence $\func{Idl}(T)$ is a chain.

In both cases, adjoining the new top $L$ gives a modular lattice (even a
distributive lattice) by Lemma \ref{mod-dist}. Thus, if $\left\vert
L\right\vert \leq 5,$ then 
\begin{equation*}
\mathfrak{K}(S_{\vee }(L))\overset{\text{Corollary \ref{ideal-join-1}}}{%
\cong }\func{Idl}(L)=\func{Idl}(T)\cup \{L\}
\end{equation*}%
is modular (even distributive) for any bounded join-semilattice $(L,\vee ,0)$
with a join-irreducible top element $1\neq 0.$ Example \ref{N5+1} supplies a
six-element join-semilattice $L:=N_{5}\cup \{\mathbf{1}\}$ with
join-irreducible $\mathbf{1}\neq $ $0$ for which $\mathfrak{K}(S_{\vee }(L))$
is non-modular and which is with \emph{minimal size} in this family.$%
\blacksquare $
\end{rem}

\subsection*{Simplicity of Semirings induced by Semilattices}

\bigskip

In what follows, we investigate when the demirings/hemirings/semiring
induced by semilattices are ideal-simple, $k$-ideal-simple and
congruence-simple.

\bigskip

\begin{punto}
Let $(L,\vee )$ be a join-semilattice. We say that $\equiv $ is a \textbf{%
join-semilattice-congruence }on $L$ iff $\equiv $ is an equivalence relation
compatible with $\vee ,$ i.e. $a\equiv b\Longrightarrow (a\vee x)\equiv
(b\vee x)$ for all $x\in L.$ With $\func{Cong}(L)$ we denote the collection
of all join-semilattice-congruences on $L.$ We say that $L$ is a \textbf{%
congruence-simple join-semilattice} iff $\func{Cong}(L)=\{\Delta _{L},\nabla
_{L}\}.$
\end{punto}

\begin{rem}
\label{chain-cc}If $L=C_{n}:a_{1}<a_{2}<\cdots <a_{n}$ is the finite chain
with $n$-elements, then $\left\vert \func{Cong}(L)\right\vert =2^{n-1}.$ So, 
$C_{2}$ is the only congruence-simple chain.
\end{rem}

\begin{thm}
\label{join-congr}Let $(L,\vee )$ be a join-semilattice.

\begin{enumerate}
\item Fix $c\in L.$ Then $\func{Cong}(D_{\vee }^{c}(L))=\func{Cong}(L).$

\item If $(L,\vee ,0)$ is a bounded join-semilattice, then $\func{Cong}%
(H_{\vee }(L))=\func{Cong}(L).$

\item If $(L,\vee ,0)$ is a bounded join-semilattice with a join-irreducible
top $1\neq 0,$ then%
\begin{equation*}
\func{Cong}(S_{\vee }(L))=\{\theta \in \func{Cong}(L)\mid \lbrack 1]_{\theta
}=\{1\}\}\cup \{\nabla _{L}\}.
\end{equation*}
\end{enumerate}
\end{thm}

\begin{Beweis}
Let $(L,\vee )$ be a join-semilattice.

\begin{enumerate}
\item Fix $c\in L.$ Compatibility with addition of $D_{\vee }^{c}(L)$ is
precisely compatibility with the underlying join of $L.$ Compatibility with
constant multiplication is automatic.

\item Choose $c=0$ and apply (1).

\item Let $(L,\vee ,0)$ be a bounded join-semilattice with a
join-irreducible top $1\neq 0.$

Let $\equiv $ $\in \func{Cong}(L)$ be such that $[1]_{\equiv }=\{1\}.$
Assume that $x\equiv y$ for $x,y\in L$ and let $z\in L.$ If $z=1$, then%
\begin{equation*}
x\cdot _{\vee }z=x\cdot _{\vee }1\equiv x\equiv y=y\cdot _{\vee }1=y\cdot
_{\vee }z.
\end{equation*}%
If $z\neq 1$, then we have the following cases:
\end{enumerate}

\begin{itemize}
\item $x=y=1.$ In this case, $x\cdot _{\vee }z=z=y\cdot _{\vee }z.$

\item $x\neq 1,$ $y\neq 1.$ In this case, $x\cdot _{\vee }z=0=y\cdot _{\vee
}z.$
\end{itemize}

The remaining cases (in which exactly one of $x,y$ is $1$)\ would contradict
the assumption $[1]_{\equiv }=\{1\}.$ Hence $\equiv $ $\in \func{Cong}%
(S_{\vee }(L)).$

Conversely, let $\theta \in \func{Cong}(S_{\vee }(L)).$ Suppose $1\theta a$
for some $a\neq 1.$ We have the following cases:

\begin{itemize}
\item $a=0.$ In this case, $(b\cdot _{\vee }1)\theta (b\cdot _{\vee }0)$
whence $b\theta 0$ for all $b\in L,$ i.e. $\theta =L\times L.$

\item $0<a<1.$ In this case, multiplying $1\theta a$ by $a$ yields $(a\cdot
_{\vee }1)\theta (a\cdot _{\vee }a),$ i.e. $a\theta 0$ and so $1\theta 0.$
Again, as in the previous case, this forces $\theta =L\times L.$
\end{itemize}

Therefore, every non-universal semiring congruence has $[1]_{\theta
}=\{1\}.\blacksquare $
\end{Beweis}

\begin{cor}
\label{join-unit-cog}Let $(L,\vee ,0)$ be a bounded join-semilattice, $%
L^{\sharp }:=L\sqcup \{1^{\sharp }\}$ where $1^{\sharp }$ is a fresh
greatest element $1^{\sharp },$ and set $\theta ^{\sharp }:=\theta \cup
\{(1^{\sharp },1^{\sharp })\}$ for any $\theta \in \func{Cong}(L).$ Then%
\begin{equation}
\func{Cong}(S_{\vee }(L^{\sharp }))=\{\theta ^{\sharp }:\theta \in \func{Cong%
}(L)\}\cup \{L^{\sharp }\times L^{\sharp }\}.  \label{join-unitized}
\end{equation}
\end{cor}

What Corollary \ref{join-unit-cog} says is that $\func{Cong}(L)$ \emph{%
embeds naturally} in $\func{Cong}(S_{\vee }(L^{\sharp }))$ as the
congruences that keep the adjoined top $1^{\sharp }$ separate, and $\func{%
Cong}(S_{\vee }(L^{\sharp }))$ is obtained by adjoining one new greatest
congruence to that embedded copy.

\begin{thm}
\label{join-simplicity}Let $(L,\vee )$ be a join-semilattice and fix $c\in
L. $

\begin{enumerate}
\item $D_{\vee }^{c}(L)$ is congruence-simple if and only if $(L,\vee )$ is
congruence-simple.

\item $D_{\vee }^{c}(L)$ is ideal-simple if and only if $\left\vert
L\right\vert =2.$

\item $D_{\vee }^{c}(L)$ is $k$-ideal-simple if and only if%
\begin{equation*}
\left\vert U_{c}\right\vert =1\quad \text{or}\quad \left\vert
U_{c}\right\vert =2\text{ and }c\text{ is \emph{minimal} in }L.
\end{equation*}
\end{enumerate}
\end{thm}

\begin{Beweis}
Let $(L,\vee )$ be a join-semilattice. and fix $c\in L.$

\begin{enumerate}
\item This follows from the equality $\func{Cong}(D_{\vee }^{c}(L))=\func{%
Cong}(L)$ in Theorem \ref{join-congr}.

\item $(\Longrightarrow )$ Assume that $D_{\vee }^{c}(L)$ is ideal-simple.
Suppose $\left\vert L\right\vert \geq 3$ and choose $x\neq c$. If $c\vee
x\in \{c,x\}$, then $\{c,x\}$ is a proper ideal. Otherwise, $\{c,c\vee x\}$
is a proper ideal. In either case, there is a proper ideal of size at least
two, so $D_{\vee }^{c}(L)$ is \emph{not} ideal-simple, a contradiction.

$(\Longleftarrow )$ Assume that $\left\vert L\right\vert =2.$ The only
ideals of $D_{\vee }^{c}(L)$ are $\{c\}$ and $L,$ whence the demiring is
ideal-simple.

\item Apply Theorem~\ref{ideals-join}. If $\left\vert U_{c}\right\vert =1,$
then the order-ideal lattice of $U_{c}$ has only one element, whence $%
D_{\vee }^{c}(L)$ is the only $k$-ideal. If $\left\vert U_{c}\right\vert =2,$
that lattice has just two elements, and the only proper subtractive ideal of 
$D_{\vee }^{c}(L)$ is $(\downarrow c]$ which has size one exactly when $c$
is minimal in $L.$ If $\left\vert U_{c}\right\vert \geq 3,$ choose $u\in
U_{c}$ with $c<u$ and $u$ not top in $U_{c}.$ Then $(\downarrow u]$ is a
proper subtractive ideal containing $\{c,u\}$, whence $\left\vert
(\downarrow u]\right\vert \geq 2$ contradicting the $k$-ideal-simplicity of $%
D_{\vee }^{c}(L).$ The meet assertions follow by duality.$\blacksquare $
\end{enumerate}
\end{Beweis}

\begin{lem}
\label{cc-k-simple}Let $(S,+,\cdot )$ be a demiring. If $S$ is
congruence-simple, then $S$ is $k$-ideal-simple.
\end{lem}

\begin{Beweis}
Let $S$ be a congruence-simple demiring. Let $I$ be any proper $k$-ideal of $%
S$ with $\left\vert I\right\vert \geq 2.$ The Bourne relation defined by $%
x\equiv _{I}y\iff x+a=y+b$ for some $a,b\in I$ forms a valid semiring
congruence. Since $S$ is congruence-simple, the congruence relation $\equiv
_{I}$ must be either the \emph{trivial} congruence $\Delta _{S}$ or the 
\emph{universal} congruence $\nabla _{S}$.

\textbf{Case I: }$\equiv _{I}=\Delta _{S}.$ Pick $a,b\in I$ with $a\neq b.$
Clearly, $a\equiv _{I}b$, a contradiction.

\textbf{Case II. }$\equiv _{I}=\nabla _{S}$. Pick any $x\in I.$ For any $%
y\in S,$ we have $x\equiv _{I}y,$ i.e. there exist $a,$ $b\in I$ with $%
y+a=x+b\in I.$ Since $I$ is a $k$-ideal, we have $y\in I.$ Thus, $%
I=S.\blacksquare $
\end{Beweis}

\begin{rem}
\label{Uc-3}Let $(L,\vee )$ be a join-semilattice. If $c\in L$ is such that $%
\left\vert U_{c}\right\vert \geq 3,$ then $U_{c}$ contains a three-element
chain $C_{3}.$ Pick any two \emph{distinct} elements $x,y\in U_{c}$ such
that $c<x,y.$ If $x,y$ are comparable, say $x<y,$ then $c<x<y.$ If $%
x\parallel y,$ then $c<x<x\vee y.\blacksquare $
\end{rem}

\begin{thm}
\label{H-CC-primes}Let $(L,\vee ,0)$ be a bounded join-semilattice. Then 
\begin{equation}
\func{Cong}(H_{\vee }(L))=\func{Cong}((L,\vee ,0))  \label{CH=CL-join}
\end{equation}%
and the following are equivalent:

\begin{enumerate}
\item $H_{\vee }(L)$ is ideal-simple;

\item $H_{\vee }(L)$ is congruence-simple;

\item $H_{\vee }(L)$ is $k$-ideal-simple;

\item $L\cong C_{2};$

\item $H_{\vee }(L)\cong \mathfrak{B}.$
\end{enumerate}
\end{thm}

\begin{Beweis}
Since multiplication is constant, compatibility with multiplication is
automatic for every congruence compatible with addition. This proves the
congruence-lattice assertions. Let $(L,\vee ,0)$ be a bounded
join-semilattice.

$(2\Longrightarrow 3)$ Any congruence-simple hemiring is $k$-ideal-simple\
(cf. Lemma \ref{cc-k-simple}).

$(3\Longrightarrow 4)$ Suppose $\left\vert L\right\vert \geq 3.$ By Remark %
\ref{Uc-3}, $L$ contains a three-element chain $C_{3}:\ 0<a<b.$ The ideal $%
(\downarrow a]$ is clearly a non-zero proper subtractive ideal of $H_{\vee
}(L),$ a contradiction. So, $L\cong C_{2}.$

$(4\Longrightarrow 5)$ By the definition of $H_{\vee }(C_{2}),$ its addition
and multiplication tables are and so $H_{\vee }(C_{2})\cong \mathfrak{B}.$

$(5\Longrightarrow 2)$ and $(5\Longrightarrow 1\Longrightarrow 3)$ are
obvious.$\blacksquare $
\end{Beweis}

The following example illustrates that (\ref{CH=CL-join}) does not extend to
semirings induced by semilattices.

\begin{ex}
\label{three-chain}Consider the chain $C_{3}:x<y<z.$ Then 
\begin{equation*}
\func{Cong}(C_{3})=\{\Delta ,\{\{x,y\},\{z\}\},\{\{x\},\{y,z\}\},\nabla \}.
\end{equation*}%
The only other equivalence relation $\tau :=\{\{x,z\},\{y\}\}$ is \emph{not}
a congruence since $x\sim z\Longrightarrow x\vee y\sim z\vee y$ whence $%
y\sim z.$ On the hand, $x\sim z\Longrightarrow x\wedge y\sim z\wedge y$
whence $x\sim y$. Either relation, together with $x\sim z$, forces $\tau
=\nabla .$

For the join construction, let $L:0<a<1.$ Its four semilattice congruences
are 
\begin{equation*}
id_{L},\quad \theta _{0a}:=\{\{0,a\},\{1\}\},\quad \theta
_{a1}:=\{\{0\},\{a,1\}\},\quad \nabla _{L}.
\end{equation*}

For $S_{\vee }(L)$, the addition and multiplication tables are 
\begin{equation*}
\begin{array}{c|ccc}
\vee & 0 & a & 1 \\ \hline
0 & 0 & a & 1 \\ 
a & a & a & 1 \\ 
1 & 1 & 1 & 1%
\end{array}%
\qquad 
\begin{array}{c|ccc}
\cdot _{\vee } & 0 & a & 1 \\ \hline
0 & 0 & 0 & 0 \\ 
a & 0 & 0 & a \\ 
1 & 0 & a & 1%
\end{array}%
.
\end{equation*}%
The congruence $\theta _{0a}$ is compatible because its non-singleton class $%
\{0,a\}$ does not contain the exceptional element $1$, and the product of
any two elements below $1$ is $0$. The other nontrivial semilattice
congruence is \emph{not} compatible: although $a\theta _{a1}1$, we have%
\begin{equation*}
a\cdot _{\vee }a=0\qquad \text{whereas}\qquad 1\cdot _{\vee }a=a,
\end{equation*}%
and $(0,a)\notin \theta _{a1}$. Therefore, $\func{Cong}(S_{\vee
}(L))=\{id_{L},\theta _{0a},\nabla _{L}\}$ and so $\left\vert \func{Cong}%
(S_{\vee }(L))\right\vert =3<4=|\func{Cong}(L)|.$

For the meet construction, let $M:$ $\perp <a<\top $, and use $0:=\top $ and 
$1:=\perp $. Its meet-semilattice congruences are $id_{M}$, $\theta _{\perp
a}=\{\{\perp ,a\},\{\top \}\}$, $\theta _{a\top }=\{\{\perp \},\{a,\top \}\}$%
, and $\nabla _{M}$. For $S_{\wedge }(M)$, the addition and multiplication
tables are 
\begin{equation*}
\begin{array}{c|ccc}
\wedge & \perp & a & \top \\ \hline
\perp & \perp & \perp & \perp \\ 
a & \perp & a & a \\ 
\top & \perp & a & \top%
\end{array}%
\qquad 
\begin{array}{c|ccc}
\cdot _{\wedge } & \perp & a & \top \\ \hline
\perp & \perp & a & \top \\ 
a & a & \top & \top \\ 
\top & \top & \top & \top%
\end{array}%
.
\end{equation*}%
Here $\theta _{a\top }$ is compatible because the multiplicative identity $%
\perp $ is a singleton class. The congruence $\theta _{\perp a}$ fails
compatibility: 
\begin{equation*}
\perp \cdot _{\wedge }a=a\qquad \text{whereas}\qquad a\cdot _{\wedge }a=\top
,
\end{equation*}%
and $(a,\top )\notin \theta _{\bot a}$. Hence $\func{Cong}(S_{\wedge
}(M))=\{id_{M},\theta _{a\top },\nabla _{M}\}$ and so $\left\vert \func{Cong}%
(S_{\wedge }(M))\right\vert =3<4=|\func{Cong}(M)|.$

In both cases, the semiring congruences are the semilattice congruences
isolating the multiplicative identity, together with the universal
congruence.$\blacksquare $
\end{ex}

\begin{theorem}
\label{join-c-simple}Let $(L,\vee ,0)$ be a bounded join-semilattice with a
join-irreducible top element $1\neq 0.$ The following are equivalent:

\begin{enumerate}
\item $S_{\vee }(L)$ is ideal-simple;

\item $S_{\vee }(L)$ is congruence-simple;

\item $S_{\vee }(L)$ is $k$-ideal-simple;

\item $L\cong C_{2};$

\item $S_{\vee }(L)\cong \mathbb{B}.$
\end{enumerate}
\end{theorem}

\begin{Beweis}
$(2\Longrightarrow 3)$ This is a well-known fact for congruence-simple
semirings (cf. \cite[Proposition 4.4]{KNT2011}).

$(3\Longrightarrow 4)\ $Assume that $S_{\vee }(L)$ is $k$-ideal-simple.
Suppose that $\left\vert L\right\vert \geq 3$ and pick $a\in L$ such that $%
0<a<1.$ The principal ideal $(\downarrow a]$ is a non-zero proper
subtractive ideal by Corollary~\ref{ideals-join-S}, a contradiction. So, $%
L=\{0,1\}\cong C_{2}.$

$(4\Longrightarrow 5)$ By definition, the addition and multiplications
tables of $S_{\vee }(C_{2})$ are%
\begin{equation*}
\begin{tabular}{ccccc}
\begin{tabular}{ccc}
$\vee $ & $0$ & $1$ \\ 
$0$ & $0$ & $1$ \\ 
$1$ & $1$ & $1$%
\end{tabular}
&  &  &  & 
\begin{tabular}{ccc}
$\cdot _{\vee }$ & $0$ & $1$ \\ 
$0$ & $0$ & $0$ \\ 
$1$ & $0$ & $1$%
\end{tabular}%
\end{tabular}%
\end{equation*}%
and so $S_{\vee }(C_{2})\cong \mathbb{B}.$

$(5\Longrightarrow 2)$ and $(5\Longrightarrow 1\Longrightarrow 3)$ are
obvious.$\blacksquare $
\end{Beweis}

\bigskip

The following example shows that the equivalence of congruence-simplicity,
ideal-simplicity and $k$-ideal-simplicity \ref{join-c-simple} of semirings
induced by bounded join-semilattices with join-irreducible tops does not
extend to demirings induced by arbitrary join-semilattices.

\begin{ex}
\label{simple-3}Consider the chain $C_{3}:x<y<z$ and choose $c=z$.

\begin{enumerate}
\item $S_{\vee }^{z}(C_{3})$ is \emph{not} congruence-simple. In fact, its
four congruences form the Boolean lattice $B_{2}:$%
\begin{equation*}
\Delta ,\quad \{\{x,y\},\{z\}\},\quad \{\{x\},\{y,z\}\},\quad \nabla .
\end{equation*}

\item Every subset containing $z$ is an ideal, i.e.%
\begin{equation*}
\mathfrak{I}(S_{\vee }^{z}(C_{3}))=\{\{z\},\{x,z\},\{y,z\},\{x,y,z\}\}.
\end{equation*}%
Since $\left\vert \{x,z\}\right\vert =2,$ it follows that $S_{\vee
}^{z}(C_{3})$ is \emph{not} ideal-simple.

\item Because $z$ is top, each of these ideals has subtractive closure $%
C_{3},$ whence $\mathfrak{K}(S_{\vee }^{z}(C_{3}))=\{C_{3}\}$ and $S_{\vee
}^{z}(C_{3})$ is $k$-ideal-simple.$\blacksquare $
\end{enumerate}
\end{ex}

\section{Kuratowski and Nucleus Semirings Induced by Semilattices}

\bigskip

Inspired by the notions of a \emph{Kuratowski map} and a \emph{nucleus map}
introduced in \cite{Nas}, we present the following two classes of semirings.

\begin{defn}
We say that a demiring $S$ is

\textbf{Kuratowski} iff $\overline{I+J}=\overline{I}+\overline{J}$ for any
ideals $I,J$ of $S;$

\textbf{nucleus} iff $\overline{I\cap J}=\overline{I}\cap \overline{J}$ for
any ideals $I,J$ of $S.$
\end{defn}

\begin{rems}
Let $(S,+,\cdot )$ be a demiring.

\begin{enumerate}
\item For any ideals $I,$ $J$ or $S,$ we have $\overline{I}+\overline{J}%
\subseteq \overline{I+J}$ and $\overline{I\cap J}\subseteq \overline{I}\cap 
\overline{J}.$

\item The equality $\overline{I}+\overline{J}=\overline{I+J}$ may fail
already for $I=S=J$ (in case $S$ is \emph{not} a pre-hemiring). Consider the
demiring $S=2\mathbb{N}=\{2,4,6,\cdots \}$ with the usual addition and
multiplication. Then%
\begin{equation*}
\overline{S}+\overline{S}=S+S=\{4,6,\cdots \}\varsubsetneqq S=\overline{%
\{4,6,\cdots \}}=\overline{S+S}.
\end{equation*}
\end{enumerate}
\end{rems}

\begin{thm}
\label{j-K-dist}Let $(L,\vee )$ be a join-semilattice.

\begin{enumerate}
\item Fix $c\in L.$ If $L$ is distributive, then $D_{\vee }^{c}(L)$ is a
Kuratowski demiring.

\item If $(L,\vee ,0)$ is bounded, then%
\begin{equation*}
H_{\vee }(L)\text{ is Kuratowski hemiring }\Longleftrightarrow \text{ }L%
\text{ is a distributive join-semilattice.}
\end{equation*}

\item If $(L,\vee ,0)$ is bounded with a join-irreducible greatest element $%
1\neq 0,$ then the following are equivalent:

\begin{enumerate}
\item $H_{\vee }(L)$ is Kuratowski hemiring;

\item $S_{\vee }(L)$ is Kuratowski semiring;

\item $L$ is a distributive join-semilattice.
\end{enumerate}
\end{enumerate}
\end{thm}

\begin{Beweis}
Let $I,J$ be ideals of $S_{\vee }(L).$ We have%
\begin{equation*}
\overline{I+J}=\{z\in L:z\leq i\vee j\text{ for some }i\in I,\ j\in J\},
\end{equation*}%
whereas 
\begin{equation*}
\overline{I}+\overline{J}=\{x\vee y:x\leq i,\ y\leq j\text{ for some }i\in
I,\ j\in J\}.
\end{equation*}%
The inclusion $\overline{I}+\overline{J}\subseteq \overline{I+J}$ always
holds.

\begin{enumerate}
\item Fix $c\in L.\ $Assume that $L$ is a distributive join-semilattice. If $%
z\in \overline{I+J}$, then $z\leq i\vee j$ for some $i\in I$ and $j\in J$.
Since $L$ is distributive, there exist $x\leq i$ and $y\leq j$ such that $%
z=x\vee y.$ By Proposition \ref{ideals-join} (2), we have $x\in \overline{I}$
and $y\in \overline{J}$ whence $z=x\vee y\in \overline{I}+\overline{J}.$ So, 
$D_{\vee }^{c}(L)$ is a Kuratowski demiring.

\item Let $(L,\vee ,0)$ be a bounded join-semilattice. Assume that $H_{\vee
}(L)$ is a Kuratowski hemiring. Let $a,b,z\in L$ be such that $z\leq a\vee
b. $ If $z\leq a$ (or $z\leq b$), then write $z=z\vee 0$ ($z=0\vee z$). So,
assume that \emph{neither} $z\leq a$ \emph{nor} $z\leq b$ (so, neither $a$
nor $b$ is top in $L$). Consider the ideals $I=\{0,a\},$ $J=\{0,b\}$ and
notice that $I+J=\{0,a,b,a\vee b\}.$ Since $H_{\vee }(L)$ is a Kuratowski
hemiring, we have%
\begin{equation*}
\overline{I+J}=\overline{I}+\overline{J}=(\downarrow I]+(\downarrow J]
\end{equation*}%
and so $z=x\vee y$ for some $x\leq a$ and $y\leq b$. Therefore $L$ is
distributive.

\item The same argument in (2) can be applied to the semiring case to show
that: $L$ is distributive $\Longleftrightarrow $ $S_{\vee }(L)$ is a
Kuratowski semiring.$\blacksquare $
\end{enumerate}
\end{Beweis}

\begin{cor}
\label{m-lattcice-dist}Let $(L,\wedge ,\vee )$ be a bounded lattice.

\begin{enumerate}
\item Fix $c\in L.$ If $L$ is distributive, then $D_{\vee }^{c}(L)$ and $%
D_{\wedge }^{c}(L)$ are Kuratowski demirings;

\item Assume that $(L,\wedge ,\vee ,0,1)$ is a bounded lattice. Then the
following are equivalent:

\begin{enumerate}
\item $H_{\vee }(L)$ is a Kuratowski hemiring;

\item $H_{\wedge }(L)$ is a Kuratowski hemiring;

\item $L$ is a distributive lattice.
\end{enumerate}

\item Assume that $(L,\wedge ,\vee ,0,1)$ is a bounded lattice with $1$
join-irreducible and $0$ meet-irreducible. Then the following are equivalent:

\begin{enumerate}
\item $S_{\vee }(L)$ is a Kuratowski semiring;

\item $S_{\wedge }(L)$ is a Kuratowski semiring;

\item $L$ is a distributive lattice.
\end{enumerate}
\end{enumerate}
\end{cor}

\begin{Beweis}
The proof follows from Theorem \ref{j-K-dist}, the isomorphisms%
\begin{equation*}
D_{\wedge }^{c}(L)\cong D_{\wedge ^{\mathrm{op}}}^{c}(L^{\mathrm{op}}),\quad
H_{\wedge }(L)\cong H_{\wedge ^{\mathrm{op}}}(L^{\mathrm{op}}),\quad
S_{\wedge }(L)\cong S_{\wedge ^{\mathrm{op}}}(L^{\mathrm{op}}),
\end{equation*}%
from Proposition \ref{semiring-meet} and the fact that a lattice $L$ is
distributive if and only if its dual lattice $L^{\mathrm{op}}$ is
distributive \cite[Corollary 62]{Gra2010}.$\blacksquare $
\end{Beweis}

The following example shows the converse of Theorem \ref{j-K-dist} (1) does 
\emph{not} hold for arbitrary $c\in L:$

\begin{ex}
\label{N5-not-Kurat}Let $L=N_{5},$ choose $c=t$ (the top of $N_{5}$) and
consider the demiring $D_{\vee }^{t}(N_{5}).$ Every ideal of $D_{\vee
}^{t}(N_{5})$ contains $t$ and so the subtractive closure of any ideal is $%
L. $ If $I,J$ are ideals of $D_{\vee }^{t}(N_{5}),$ then $\overline{I+J}%
=L=L+L=\overline{I}+\overline{J}.$ So, $D_{\vee }^{t}(N_{5})$ is a
Kuratowski demiring. However, $N_{5}$ is non-distributive (in fact
non-modular).$\blacksquare $
\end{ex}

While a sublattice of modular (distributive) lattice is also modular
(distributive). These facts does not extend to join-semilattices.

\begin{ex}
Consider the eight-element Boolean lattice $B_{3}=\mathcal{P}(\{1,2,3\})$
ordered by inclusion, with join given by union and meet by intersection.
Then $B_{3}$ is a distributive (and modular) lattice. Set 
\begin{equation*}
p=\{1,2\},\qquad q=\{1,3\},\qquad r=\{2,3\},\qquad \top =\{1,2,3\}.
\end{equation*}%
Then $L:=\{p,q,r,\top \}\subseteq B_{3}$ is closed under unions, i.e. $L$ is
a sub-join-semilattice of the distributive lattice $B_{3}$.

\textbf{Claim I:}\ $L$ is \emph{non-modular}. Take $a=b=p$ and $z=q$. Then $%
p\leq p\leq p\vee q=\top .$ But the only element of $L$ below $q$ is $q$,
and $p\vee q=\top \neq p$. Thus no modularity witness exists in $L$.

\textbf{Claim II:\ }$L$ is \emph{non-distributive}. Notice that $p\leq q\vee
r=\top .$ The only element of $L$ below $q$ is $q$, and the only element of $%
L$ below $r$ is $r$. Consequently, there are no $u,v\in L$ with $u\leq q$, $%
v\leq r$, and $p=u\vee v$, since $q\vee r=\top \neq p.\blacksquare $
\end{ex}

Although the modularity (distributivity)\ of a join-semilattice is not
necessarily inherited by arbitrary sub-join-semilattices, the upper cone $%
U_{c}$ is a special case.

\begin{prop}
\label{Uc-hereditary}Let $(L,\vee )$ be a join-semilattice and fix $c\in L.$

\begin{enumerate}
\item If $L$ is a modular join-semilattice, then $U_{c}$ is modular.

\item If $L$ is a distributive join-semilattice, then $U_{c}$ is
distributive.
\end{enumerate}
\end{prop}

\begin{Beweis}
Consider the upper cone $U_{c}:=[c\uparrow ).$

\begin{enumerate}
\item Let $a,b,z\in U_{c}$ satisfy $a\leq b\leq a\vee z.$ By modularity of $%
L,$ there exists $d\in L$ such that $d\leq z$ and $b=a\vee d.$ If $d\in
U_{c},$ we are done. If $d\notin U_{c},$ set $d^{\prime }:=c\vee d\in U_{c}.$
Because $d\leq z$ and $c\leq z$, we have $d^{\prime }\leq z$. Since $c\leq a$%
, we have%
\begin{equation*}
a\vee d^{\prime }=a\vee (c\vee d)=(a\vee c)\vee d=a\vee d=b.
\end{equation*}%
Hence $d^{\prime }$ is a witness to modularity within $U_{c}.$

\item Let $x,y,z\in U_{c}$ satisfy $x\leq y\vee z.$ By distributivity of $L$%
, there exist $u,v\in L$ such that 
\begin{equation*}
u\leq y,\qquad v\leq z,\qquad x=u\vee v.
\end{equation*}%
If $u,v\in U_{c},$ we are done. If not, set $u^{\prime }:=u\vee c\in U_{c}$
and $v^{\prime }:=v\vee c\in U_{c}.$ Since $c\leq y$ and $c\leq z$, we have $%
u^{\prime }\leq y$ and $v^{\prime }\leq z.$ Finally, since $c\leq x$, we have%
\begin{equation*}
u^{\prime }\vee v^{\prime }=(u\vee c)\vee (v\vee c)=(u\vee v)\vee c=x\vee
c=x.
\end{equation*}%
Thus the required decomposition takes place inside $U_{c}$.$\blacksquare $
\end{enumerate}
\end{Beweis}

The reader should be warend that, despite the fact that modularity and
distributivity\ are herditary properties for lattices (which follows
directly from Lemma \ref{mod-dist}), none of them is so for
join-semilattices. The fact that $U_{c}$ inherits the modulariy
(distributivity)\ of $L,$ as shown above, stems from the structure of $U_{c}$
as a prinicpal up-set and does not generalize to arbitrary
sub-join-semilattices of join-lattices.

\begin{ex}
Consider the eight-element Boolean lattice $B_{3}:=\mathcal{P}(\{1,2,3\})$
ordered by inclusion and the sub-join-semilattice $T=\{\{1,2\},\{1,3\},\{2,3%
\},\{1,2,3\}\}.$ Clearly, $B_{3}$ is distributive (whence modular), while $T$
is not even modular:\ let $u=w=\{1,2\}$ and $v=\{1,3\}.$ Then $u\leq w\leq
u\vee v,$ but thre exists \emph{no} $y\in T$ such that $y\leq v$ and $%
w=u\vee y.\blacksquare $
\end{ex}

The following example demonstrates that the distributivity of $U_{c}$ alone
is \emph{not sufficient} to guarantee that $D_{\vee }^{c}(L)$ is a
Kuratowski demiring in Theorem \ref{j-K-dist} (1):

\begin{ex}
\label{Uc-not-suff}Consider the three-element join-semilattice 
\begin{equation*}
V:p\parallel q,\text{ }p\vee q=t,
\end{equation*}%
and choose $c=p.$ Its join and multiplication tables are 
\begin{equation*}
\begin{array}{c|ccc}
\vee & p & q & t \\ \hline
p & p & t & t \\ 
q & t & q & t \\ 
t & t & t & t%
\end{array}%
\qquad 
\begin{array}{c|ccc}
\cdot _{p} & p & q & t \\ \hline
p & p & p & p \\ 
q & p & p & p \\ 
t & p & p & p%
\end{array}%
.
\end{equation*}%
The upper cone $U_{p}=\{p,t\}\cong C_{2}$ is a chain whence distributive,
but $D_{\vee }^{p}(V)$ is \emph{not} a Kuratowski demiring. Let $I=\{p\}$
and $J=\{p,t\}$ and notice that $I+J=J.$ So, we have%
\begin{equation*}
\overline{I}+\overline{J}=I+V=\{p,t\}\varsubsetneqq V=\overline{J}=\overline{%
I+J}.\blacksquare
\end{equation*}
\end{ex}

\begin{thm}
\label{nuc-dem}Let $(L,\vee )\ $be a join-semilattice.

\begin{enumerate}
\item Fix $c\in L.$ Then $D_{\vee }^{c}(L)$ is a nucleus demiring if and
only if $\left\vert U_{c}\right\vert \leq 2.$

\item Assume that $(L,\vee ,0)$ is a bounded join-semilattice. Then $H_{\vee
}(L)$ is a nucleus hemiring if and only if $\left\vert L\right\vert =2.$
\end{enumerate}
\end{thm}

\begin{Beweis}
Let $(L,\vee )\ $be a join-semilattice.

\begin{enumerate}
\item Fix $c\in L.$ $(\Longrightarrow )$ Assume that $D_{\vee }^{c}(L)$ is a
nucleus demiring and suppose that $\left\vert U_{c}\right\vert \geq 3.$ By
Remark \ref{Uc-3}, $U_{c}$ contains a three-element chain $C_{3}:$ $c<a<b.$
Set $I:=\{c,a\}$ and $J:=\{c,b\}$. Then $a\in \overline{I}\cap \overline{J}$
while $a\notin \overline{I\cap J}=\overline{\{c\}}=(\downarrow c].$ So, $%
\overline{I\cap J}\varsubsetneqq \overline{I}\cap \overline{J},$ a
contradiction.

$(\Longleftarrow )$ Assume that $\left\vert U_{c}\right\vert \leq 2.$

\textbf{Case I:}\ $U_{c}=\{c\}.$ In this case, $c$ is the top of $L,$ and $%
\overline{I}=L$ for every ideal $I$ of $D_{\vee }^{c}(L).$ Clearly, $D_{\vee
}^{c}(L)$ is nucleus in this case.

\textbf{Case II:}$\ U_{c}=\{c,u\}$ with $c<u.$ Then $u$ is the top of $L.$
For any ideal $K$ of $D_{\vee }^{c}(L),$ \emph{either} $K\cap U_{c}=\{c\}$
whence $\overline{K}=(\downarrow c],$ \emph{or} $I\cap U_{c}=\{c,u\}$ whence 
$\overline{I}=L.$

Let $I$ and $J$ be ideals of $D_{\vee }^{c}(L).$

\begin{itemize}
\item $I\cap U_{c}=\{c,u\}=J\cap U_{c}.$ In this case, we have $\overline{%
I\cap J}=L=L\cap L=\overline{I}\cap \overline{J}.$

\item $I\cap U_{c}=\{c,u\}$ and $J\cap U_{c}=\{c\}.$ In this case, $%
\overline{I\cap J}=(\downarrow c]=L\cap (\downarrow c]=\overline{I}\cap 
\overline{J}.$

\item $I\cap U_{c}=\{c\}$ and $J\cap U_{c}=\{c,u\}.$ In this case, $%
\overline{I\cap J}=(\downarrow c]=(\downarrow c]\cap L=\overline{I}\cap 
\overline{J}.$

\item $I\cap U_{c}=\{c\}=J\cap U_{c}=\{c\}.$ In this case, $\overline{I\cap J%
}=(\downarrow c]=(\downarrow c]\cap (\downarrow c]=\overline{I}\cap 
\overline{J}.$
\end{itemize}

In all subcases, $\overline{I\cap J}=\overline{I}\cap \overline{J},$ whence $%
D_{\vee }^{c}(L)$ is nucleus.

\item Choose $c=0$ and apply (1).$\blacksquare $
\end{enumerate}
\end{Beweis}

\begin{thm}
\label{j-Nuc-dist}Let $(L,\vee )$ be a bounded join-semilattice with a
join-irreducible top $1\neq 0.$ The following are equivalent:

\begin{enumerate}
\item $S_{\vee }(L)$ is a nucleus semiring;

\item $(\downarrow a]\cap (\downarrow b]=\{0\}$ for all distinct $a,b\in
L\backslash \{0,1\};$

\item $\left\vert L\right\vert =2$ or $\left\vert L\right\vert =3.$
\end{enumerate}
\end{thm}

\begin{Beweis}
$(1\Longrightarrow 2)$ Assume that $S_{\vee }(L)$ is a nucleus semiring. For
distinct $a,b\in L\backslash \{0,1\},$ we have for the ideals $I:=\{0,a\}$
and $J=\{0,b\}:$%
\begin{equation*}
(\downarrow a]\cap (\downarrow b]=\overline{I}\cap \overline{J}=\overline{%
I\cap J}=\overline{\{0\}}=\{0\}.
\end{equation*}

$(2\Longrightarrow 3)$ Suppose that $\left\vert L\right\vert \geq 4.$ Let $%
a,b\in L\backslash \{0,1\}$ be distinct elements of $L.$

\textbf{Case I:} $a$ and $b$ are comparable, say $a\leq b.$ In this case, $%
0\neq a\in (\downarrow a]\cap (\downarrow b],$ a contradiction.

\textbf{Case II:} $a\parallel b.$ In this case, $c:=a\vee b\notin \{a,1\}$
(since $1$ is join-irreducible) and $0\neq a\in (\downarrow a]\cap
(\downarrow c],$ a contradiction.

$(3\Longrightarrow 1)$ Assume $\left\vert L\right\vert \leq 3.$ Then $L$ is
a chain:

\textbf{Case I:\ }$L=C_{2}:\ 0<1.$ In this case, $S_{\vee }(L)\cong \mathbb{B%
}$ which is clearly nucleus.

\textbf{Case II:\ }$L=C_{3}:$ $0<a<1$ for some $a\in L.$ In this case, $%
\mathfrak{I}(S_{\vee }(L))=\{\{0\},\{0,a\},L\}=\mathfrak{K}(S_{\vee }(L))$,
i.e. $S_{\vee }(L)$ is subtractive whence nucleus.$\blacksquare $
\end{Beweis}

It was claimed in \cite[Lemma 2.1 (7)]{GD2024} that the subtractive closure
of an (arbitrary) intersection of ideals of any commutative semiring is
equal to the intersection of the subtractive closures of these ideals.
Theorem \ref{j-Nuc-dist} provides us not only with infinitely many
commutative semirings for which this claim does not hold, but also provides
us with a machine to produce concrete counterexamples which is a clear
advantage over the recently given counterexamples \cite[Remark 2.6]{SGBS}
and \cite[Theorem 3.6]{Nas}).

\bigskip

The following example shows that in the assumption $\left\vert L\right\vert
\leq 3$ that guarantees $S_{\vee }(L)$ is nucleus, the upper bound $3$ is
sharp.

\bigskip

\begin{ex}
\label{C4-not-N}Consider the chain 4-element chain%
\begin{equation*}
C_{4}:0<a<b<1
\end{equation*}%
and the semiring $S_{\vee }(C_{4}).$ Consider the ideals $I:=\{0,a\}$ and $%
J:=\{0,b\}$ and notice that $\overline{I}=I$ while $\overline{J}=\{0,a,b\}.$
Then%
\begin{equation*}
\overline{I\cap J}=\overline{\{0\}}=\{0\}\neq \{0,a\}=\{0,a\}\cap \{0,a,b\}=%
\overline{I}\cap \overline{J}.
\end{equation*}%
So, $S_{\vee }(C_{4})$ is \emph{not} a nucleus semiring as predicted by
Theorem \ref{j-Nuc-dist}. A direct check shows that $S_{\vee }(C_{4})$ is a
Kuratowski semiring as predicted by Corollary \ref{m-lattcice-dist} (since $%
\func{Idl}(L)$ is distributive by Lemma \ref{mod-dist}).
\end{ex}

\bigskip

\begin{ex}
\label{N5+1-no-KN}Consider the semiring $S_{\vee }(N_{5}\cup \{1^{\sharp
}\}) $ from Example \ref{N5+1} and the ideals $I=\{e,c\},$ $J=\{e,a\}$ and $%
K=\{e,b\}.$ Then we have%
\begin{eqnarray*}
\overline{I\cap J} &=&\overline{\{e\}}=\{e\}\varsubsetneqq
\{e,a\}=\{e,a,c\}\cap \{e,a\}=\overline{I}\cap \overline{J}; \\
\overline{J+K} &=&\overline{\{e,a,b,t\}}=N_{5}\varsupsetneqq
\{e,a,b,t\}=\{e,a\}+\{e,b\}=\overline{J}+\overline{K}.
\end{eqnarray*}%
So, $S_{\vee }(N_{5}\cup \{1^{\sharp }\})$ are \emph{neither} a Kuratowski
semiring \emph{nor} a nucleus semiring.$\blacksquare $
\end{ex}

\section{Noetherian and Artinian properties}

\bigskip

Let $S$ be a commutative demiring. A prime ideal $P\subseteq S$ is a \emph{%
proper} ideal $P\varsubsetneqq S$ such that for any $s,t\in S$ with $st\in
P, $ we have either $s\in P$ or $t\in P.$ We denote with $\func{Spec}(S)$
(resp. $\func{Spec}_{k}(S)$)\ the class of prime ideals (resp. subtractive
prime ideals)\ of $S.$ We define the \textbf{Krull dimension} $K.\dim (S)$ (%
\textbf{Krull} $k$\textbf{-dimension} $K.\dim _{k}(S)$)\ of a demiring $S$
to be we the supremum of the lengths of chains of prime ideals (\emph{%
subtractive} prime ideals) of $S.$ Indeed $K.\dim _{k}(S)\leq K.\dim (S)$
for any commutative demiring $S.$

\bigskip

\begin{defn}
We call a commutative demiring $S$

\textbf{Noetherian} (\textbf{Artinian})\ iff the collection $\mathfrak{I}%
(S)\ $of all ideals of $S$ satisfies the ACC (DCC);

$k$\textbf{-Noetherian} ($k$\textbf{-Artinian}) iff the collection $%
\mathfrak{K}(S)\ $of all subtractive ideals of $S$ satisfies the ACC (DCC).
\end{defn}

\begin{theorem}
\label{prime-collapse}

\begin{enumerate}
\item Let $(L,\vee ,0)\ $be a bounded join-semilattice with a
join-irreducible greatest element $1\neq 0.$ Then%
\begin{equation*}
\func{Spec}(S_{\vee }(L))=\func{Spec}_{k}(S_{\vee }(L))=\{L\backslash
\{1_{L}\}\}.
\end{equation*}%
In particular, we have%
\begin{equation*}
K.\dim (S_{\vee }(L))=0=K.\dim _{k}(S_{\vee }(L)).
\end{equation*}

\item Let $(M,\wedge ,\top )\ $be a bounded meet-semilattice with a
meet-irreducible least element $\bot \neq \top .$ Then%
\begin{equation*}
\func{Spec}(S_{\wedge }(M))=\{M\backslash \{\bot \}\}=\func{Spec}%
_{k}(S_{\wedge }(M)).
\end{equation*}%
In particular, we have%
\begin{equation*}
K.\dim (S_{\wedge }(M))=0=K.\dim _{k}(S_{\wedge }(M)).
\end{equation*}
\end{enumerate}
\end{theorem}

\begin{Beweis}
\begin{enumerate}
\item Set $\mathfrak{m}:=L\backslash \{1\}.$ Clearly, $0\in \mathfrak{m}$
and the join-irreducibility of $1$ implies that $\mathfrak{m}$ is
join-closed. So, $\mathfrak{m}$ is a proper ideal of $S_{\vee }(L).$ If $%
a,b\notin \mathfrak{m}$, then $a=b=1$ and so $a\cdot _{\vee }b=1_{L}\notin 
\mathfrak{m}.$ Thus $\mathfrak{m}$ is a prime ideal of $S_{\vee }(L).$

Conversely, let $P$ be any prime ideal of $S_{\vee }(L).$ Since $P$ is a
proper ideal, $1\notin P.$ For every $x\in L\backslash \{1\},$ we have $%
x^{2}=0\in P,$ whence $x\in P.$ Hence $L\backslash \{1\}\subseteq P\subseteq
L\backslash \{1\},$ whence $P=L\backslash \{1\}.$ Clearly, $L\backslash
\{1\}\ $is a down set and closed under joins, whence $L\backslash \{1\}\in 
\func{Idl}(L)\overset{\text{Corollary \ref{ideal-join-1}}}{\cong }K(S_{\vee
}(L)).$ So, we have%
\begin{eqnarray*}
\func{Spec}(S_{\vee }(L)) &=&\{L\backslash \{1\}\}=\func{Spec}_{k}(S_{\vee
}(L)); \\
K.\dim (S_{\vee }(L)) &=&0=K.\dim _{k}(S_{\vee }(L)).
\end{eqnarray*}

\item Using the duality in Proposition \ref{semiring-meet}, we have%
\begin{equation*}
\begin{tabular}{ccccc}
$\func{Spec}(S_{\wedge }(M))$ & $=$ & $\func{Spec}(S_{\wedge ^{\mathrm{op}%
}}(M^{\mathrm{op}}))$ & $=$ & $\func{Spec}_{k}(S_{\wedge ^{\mathrm{op}}}(M^{%
\mathrm{op}}))$ \\ 
& $=$ & $\{M\backslash \{1_{M^{\mathrm{op}}}\}\}$ & $=$ & $\{M\backslash
\{\bot \}\}.$%
\end{tabular}%
\end{equation*}%
and so%
\begin{equation*}
K.\dim (S_{\wedge }(M))=0=K.\dim _{k}(S_{\wedge }(M)).\blacksquare
\end{equation*}
\end{enumerate}
\end{Beweis}

As a direct consequence of Proposition \ref{ideals-join} and Corollary \ref%
{ideals-join-S}, we have:

\begin{theorem}
\label{j-semiring-NA}Let $(L,\vee )$ be a join-semilattice.

\begin{enumerate}
\item Fix $c\in L$.

\begin{enumerate}
\item $D_{\vee }^{c}(L)$ is Noetherian (Artinian)\ if and only if $\mathfrak{%
I}_{\vee }^{c}(L)$ satisfies the ACC (DCC).

\item $D_{\vee }^{c}(L)$ is $k$-Noetherian ($k$-Artinian) if and only if $%
\func{Idl}(U_{c})$ satisfies the ACC (DCC).
\end{enumerate}

\item Assume that $(L,\vee ,0)$ is a bounded join-semilattice.

\begin{enumerate}
\item $H_{\vee }(L)$ is Noetherian (Artinian)\ if and only if $\mathfrak{I}%
_{\vee }(L):=\mathfrak{I}_{\vee }^{0}(L)$ satisfies the ACC (DCC).

\item $H_{\vee }(L)$ is $k$-Noetherian ($k$-Artinian)\ if and only if $\func{%
Idl}(L)$ satisfies the ACC (DCC).
\end{enumerate}

\item Assume that $(L,\vee ,0)$ is a bounded join-semilattice with a
join-irreducible $1\neq 0.$

\begin{enumerate}
\item $S_{\vee }(L)$ is Noetherian (Artinian)\ if and only if $\mathfrak{I}%
_{\vee }(L)$ satisfies the ACC (DCC).

\item $S_{\vee }(L)$ is $k$-Noetherian ($k$-Artinian) if and only if $\func{%
Idl}(L)$ satisfies the ACC (DCC).
\end{enumerate}
\end{enumerate}
\end{theorem}

\begin{cor}
\label{chain-4-types}Let $(C,\leq ;0_{C},1_{C})$ be a non-trivial bounded
chain and consider the associated semirings 
\begin{equation*}
S_{\vee }(C):=(C;\max ,0_{C};\cdot _{\vee },1_{C}),\qquad S_{\wedge
}(C):=(C;\min ,1_{C};\cdot _{\wedge },0_{C}).
\end{equation*}

\begin{enumerate}
\item $\mathfrak{K}(S_{\vee }(C))$ and $\mathfrak{K}(S_{\wedge }(C))$ are
distributive (and modular).

\item The following are equivalent:

\begin{enumerate}
\item $C$ is finite;

\item $S_{\vee}(C)$ is Noetherian;

\item $S_{\vee}(C)$ is Artinian;

\item $S_{\wedge}(C)$ is Noetherian;

\item $S_{\wedge}(C)$ is Artinian.
\end{enumerate}

\item $S_{\vee}(C)$ is $k$-Noetherian ($k$-Artinian) if and only if $C$
satisfies the ACC (DCC).

\item $S_{\wedge}(C)$ is $k$-Noetherian ($k$-Artinian) if and only if $C$
satisfies the DCC (ACC).
\end{enumerate}
\end{cor}

\begin{Beweis}
Since $C$ is a non-trivial chain, its greatest element $1_{C}$ is
join-irreducible and its least element $0_{C}$ is meet-irreducible: $\max
(x,y)=x\vee y=1$ forces $x=1$ or $y=1$ and $\min (x,y)=x\wedge y=0$ forces $%
x=0$ or $y=0.$ Therefore $S_{\vee }(C)$ and $S_{\wedge }(C)$ are additively
idempotent semirings.

\begin{enumerate}
\item By Corollary \ref{ideals-join-S} (2) and its dual, we have%
\begin{equation*}
\mathfrak{K}(S_{\vee }(C))\cong \func{Idl}(C),\qquad \mathfrak{K}(S_{\wedge
}(C))\cong \func{Fil}(C).
\end{equation*}%
The order ideals of the chain $C$ are its non-empty initial segments, and
the filters are its non-empty final segments. Each collection is linearly
ordered by inclusion. Hence both lattices are modular and distributive.

\item Since every subset of any chain is closed under binary maxima and
minima, Propositions \ref{ideals-join}(1) and its dual give%
\begin{eqnarray*}
\mathfrak{I}(S_{\vee }(C)) &:&=\{C\}\cup \{A\subseteq C\backslash \{1_{C}\}%
\text{ }\mid \text{ }0_{C}\in A\}; \\
\mathfrak{I}(S_{\wedge }(C)) &:&=\{C\}\cup \{B\subseteq C\backslash \{0_{C}\}%
\text{ }\mid \text{ }1_{C}\in B\}.
\end{eqnarray*}%
Indeed, $1_{C}$ is the multiplicative identity in $S_{\vee }(C)$, while $%
0_{C}$ is the multiplicative identity in $S_{\wedge }(C)$. Thus the only
ideal of $S_{\vee }(C)$ that contains $1_{C}$ is $S_{\vee }(C),$ and the
only ideal of $S_{\wedge }(C)$ that contains $0_{C}$ is $S_{\wedge }(C).$ If 
$C$ is finite, then both semirings are finite, and hence both are Noetherian
and Artinian.

Assume that $C$ is infinite and set $C^{\circ }:=C\backslash \{0,1\}.$
Choose pairwise distinct elements $c_{1},c_{2},\cdots \in C^{\circ }$ (with
no monotonicity assumption). For $n\geq 1$, define%
\begin{eqnarray*}
G_{n} &:&=\{0,c_{1},\cdots ,c_{n}\},\qquad I_{n}:=\{0,c_{n+1},c_{n+2},\cdots
\}; \\
J_{n} &:&=\{1,c_{1},\cdots ,c_{n}\},\qquad K_{n}:=\{1,c_{n+1},c_{n+2},\cdots
\}.
\end{eqnarray*}%
Each $G_{n}$ and $I_{n}$ contains $0$ and omits $1$, whence a proper ideal
of $S_{\vee }(C)$. Each $J_{n}$ and $K_{n}$ contains $1$ and omits $0$,
whence a proper ideal of $S_{\wedge }(C)$. All these sets are closed under
the relevant addition because the maximum or minimum of two elements of a
chain is one of those elements. Since the $c_{k}$ are distinct, we have the
following strict chains%
\begin{eqnarray*}
G_{1} &\varsubsetneqq &G_{2}\varsubsetneqq \cdots ,\quad I_{1}\varsupsetneqq
I_{2}\varsupsetneqq \cdots ,\quad \\
J_{1} &\varsubsetneqq &J_{2}\varsubsetneqq \cdots ,\quad K_{1}\varsupsetneqq
K_{2}\varsupsetneqq \cdots .
\end{eqnarray*}%
Thus both semirings $S_{\vee }(C)$ and $S_{\wedge }(C)$ are neither
Noetherian nor Artinian.

\item By Proposition \ref{ideals-join} (3), $\mathfrak{K}(S_{\vee }(C))\cong 
\func{Idl}(C)$, the chain of non-empty initial segments of $C$. This lattice
satisfies the ACC if and only if $C$ does: any strictly ascending sequence $%
a_{1}<a_{2}<\cdots $ in $C$ gives a strictly ascending sequence of principal
initial segments. Conversely, for every strictly ascending sequence of
initial segments $D_{1}\varsubsetneqq D_{2}\varsubsetneqq \cdots
\varsubsetneqq D_{n}\varsubsetneqq D_{n+1}\cdots $ choose $x_{n}\in
D_{n+1}\backslash D_{n}.$ For $m>n,$ the element $x_{n}\in D_{m}$ but $%
x_{m}\notin D_{m}.$ If $x_{m}\leq x_{n},$ the down-closure of $D_{m}$ would
imply $x_{m}\in D_{m},$ a contradiction. So, we have a strictly ascending
sequence in $C:x_{1}<x_{2}<\cdots <x_{n}<x_{n+1}\cdots .$ The corresponding
DCC equivalence follows dually. This proves the assertion.

\item $\mathfrak{K}(S_{\wedge }(C))\cong \mathfrak{K}(S_{\wedge ^{\mathrm{op}%
}}(C^{\mathrm{op}}))\cong \func{Idl}(C^{\mathrm{op}})\cong \func{Fil}(C)$
(the chain of non-empty final segments of $C$). This is the order-dual of $%
\func{Idl}(C)$, so it satisfies the ACC (respectively DCC) if and only if $C$
satisfies the DCC (respectively ACC).$\blacksquare $
\end{enumerate}
\end{Beweis}

\begin{exs}
\label{k-N-k-A}

\begin{enumerate}
\item Consider the bounded chain%
\begin{equation*}
\mathbb{W}\cup \{\infty \}:0<1<2<\cdots <n<n+1<\cdots <\infty .
\end{equation*}%
Since the chain is infinite, $S_{\vee }(\mathbb{W}\cup \{\infty \})$ and $%
S_{\wedge }(\mathbb{W}\cup \{\infty \})$ are \emph{neither} Noetherian \emph{%
nor} Artinian by Corollary \ref{chain-4-types} (2). Since the chain
satisfies the DCC but not the ACC, it follows by Corollary \ref%
{chain-4-types} (3) that $S_{\vee }(\mathbb{W}\cup \{\infty \})$ is $k$%
-Artinian but \emph{not }$k$-Noetherian, while $S_{\wedge }(\mathbb{W}\cup
\{\infty \})$ is $k$-Noetherian but \emph{not} $k$-Artinian.

\item Consider the bounded chain $C=\{0\}\cup \{\frac{1}{n}$ $\mid $ $\in 
\mathbb{N}\}:$ 
\begin{equation*}
C:0<\cdots <\frac{1}{n+2}<\frac{1}{n+1}<\cdots <\frac{1}{4}<\frac{1}{3}<%
\frac{1}{2}<1.
\end{equation*}%
Since $C$ is infinite, $S_{\vee }(C)$ and $S_{\wedge }(C)$ are \emph{neither}
Noetherian \emph{nor} Artinian by Corollary \ref{chain-4-types} (2). Since $%
C $ satisfies the ACC condition but not the DCC condition, it follows by
Corollary \ref{chain-4-types} (3) that $S_{\vee }(C)$ is $k$-Noetherian but 
\emph{not }$k$-Artinian, while $S_{\wedge }(C)$ is $k$-Artinian but \emph{not%
} $k$-Noetherian.

\item Consider the unit interval $[0,1]$ its usual order. Since the infinite
chain $[0,1]$ satisfies neither the ACC nor the DCC, it follows by Corollary %
\ref{chain-4-types} (3) that $S_{\vee }([0,1])$ and $S_{\wedge }([0,1])$ are 
\emph{neither} $k$-Noetherian \emph{nor} $k$-Artinian.
\end{enumerate}
\end{exs}

\subsection*{Declarations}

\medskip

\begin{itemize}
\item \textbf{Funding:\ }No funding was received to assist with the
preparation of this manuscript.

\item \textbf{Financial Interests:\ }The author has no relevant financial or
non-financial interests to disclose.

\item \textbf{Usage of AI:\ }The manuscript was audited by \emph{ChatGPT 6
Luna} and \emph{Claude Sonnet 5.5}. The author is solely responsible for the
correctness of all results and examples/counterexamples in this manuscript.
\end{itemize}

\end{document}